\documentclass [reqno,10pt,oneside]{amsart}
\usepackage{amsthm}
\usepackage{amscd}
\usepackage{amsmath,amssymb}
\usepackage{algorithm}
\usepackage{placeins}
\usepackage{float}
\usepackage[noend]{algorithmic}
\usepackage{tikz}
\usetikzlibrary{matrix}

\usepackage{color}
\usepackage[colorlinks, plainpages]{hyperref}

\theoremstyle {plain}
\newtheorem {thm}{Theorem}[section]
\newtheorem {prop}[thm]{Proposition}
\newtheorem {lem}[thm]{Lemma}
\newtheorem {cor}[thm]{Corollary}
\theoremstyle {definition}
\newtheorem {defn}[thm]{Definition}

\theoremstyle {remark}
\newtheorem {rem}[thm]{Remark}

\newtheorem {exmp}[thm]{Example}

\DeclareMathOperator{\ord}{ord}
\DeclareMathOperator{\Max}{Max}

\newcommand{\C}{{\mathbb C}}

\newcommand{\N}{{\mathbb N}}
\newcommand{\Q}{{\mathbb Q}}

\newcommand{\mA}{{\mathcal A}}

\begin{document}

\bibliographystyle{alpha}

\title{Classification of Unimodal Parametric Plane Curve Singularities in Positive Characteristic}

\author{Muhammad Ahsan Binyamin$^{\ast}$, Gert-Martin Greuel, Khawar Mehmood, Gerhard Pfister}

\address{Muhammad Ahsan Binyamin\\ Department of Mathematics, GC University, Faisalabad,  Pakistan}
\email{ahsanbanyamin@gmail.com}

\address{Khawar Mehmood\\ Department of Applied Mathematics and Statistics, Institute of Space and Technology, Islamabad, Pakistan}
\email{khawar1073@gmail.com}

\address{Gert-Martin Greuel, Gerhard Pfister\\ Department of Mathematics, University of Kaiserslautern, Germany}
\email {greuel@mathematik.uni-kl.de, pfister@mathematik.uni-kl.de}



\maketitle

\begin{abstract}
 In 2011 Hefez and Hernandes completed Zariski's analytic classification of plane branches belonging to a given equisingularity class by creating "very short" parameterizations over the complex numbers. Their results were used by Mehmood and Pfister to  classify unimodal plane branches in characteristic $0$ by giving lists of normal forms. The aim of this paper is to give a complete classification of unimodal plane branches over an algebraically closed field of positive characteristic. Since the methods of Hefez and Hernandes cannot be used in positive characteristic, we use a different approach and, for some sporadic singularities in small characteristic, computations with {\sc Singular}. 
Our methods  are characteristic-independent and provide a different proof of the classification in characteristic $0$ showing at the same time that this classification holds also in large characteristic. The main theoretical ingredients are the semicontinuity of the semigroup and of the modality, which we prove and which may be of independent interest.
\end{abstract}

\bigskip
\tableofcontents

\section{Introduction and Preliminaries}
\label{sec.intro}
The classification of singularities  up to some equivalence relation has a long history with important contributions going back to Arnold and Zariski (cf.  \cite{A}, \cite{ZO}). 
For fixed small invariants\,\footnote{\,For large invariants it is in general not possible to give explicit lists. Then one has to construct moduli spaces, see e.g \cite{GP0}.} such a classification means finding a list of germs with the given invariants and showing that all germs with these invariants are equivalent to a germ in the list. 

The idea of using the modality (see Definition \ref{def.modal}) as the main invariant for the classification of singularities is due to V. I. Arnold.  It has turned out to capture the most important structural properties of a singularity, at least for small modalities. E.g. Arnold's 
classification of simple (i.e., of modality 0) isolated hypersurface singularities  with respect to right equivalence leads to the famous list of ADE-singularities.

Two main approaches have been used for the classification of curve singularities defined over an algebraically closed field. The first approach is to use the defining equations for the classification (with respect to right or contact equivalence) with important contributions by Arnold \cite{A}, Guisti \cite{GM83}, Wall \cite{Wal83} and many others (\cite{DAF},\cite{BGM},\cite{BGM1},\cite{GK},\cite{GH}).
 The second approach is to use the parameterization of a curve singularity for the classification.   Bruce and Gaffney (\cite{BG}) resp. Gibson and Hobbs  (\cite{GH}) classified the simple plane resp. space curve singularities. Arnold \cite{A1}
 gave a classification for simple curves in $n$-space with a characterization by invariants  in
  \cite{BK}. There are several further contributions in characteristic  zero 
  (cf.\cite{BM}, \cite{HHS}, \cite{HHA}, \cite{JP}, \cite{IG}, \cite{MB}, \cite{MP}). 
 In \cite{HHA}, Hefez and Hernandes introduced a general approach using differential forms for the classification of complex analytic parametric plane curve singularities. Using values of differentials, they derived a "very short" parameterization, improving Zariski's short parameterization, but  without giving explicit lists.
In \cite{MP}  Mehmood and Pfister  used another approach to classify the simple and unimodal  parametric plane branches in characteristic $0$ and gave explicit lists (see also  \cite{BM} for a classifier). They also classified simple parametric plane curve singularities in positive characteristic, see \cite{MP1}. 
 
 In this paper
we  classify the unimodal  parametric plane curve singularities in positive characteristic. The main ingredients and tools for our proofs are 
\begin{itemize}
\item the semicontinuity of the $\delta$-invariant in families of parameterizations
(Theorem \ref{thm.semicont}),
\item a proper definition of $\mA$-modality and the proof of its semicontinuity (Theorem \ref{thm.semicontmod}),
\item the semicontinuity of the semigroup with few generators (Proposition \ref{prop.semicontG})
(the general semicontinuity is posed as an open problem in Remark \ref{rem.sem}),
\item a generalization of Zariski's short parameterization to arbitrary characteristic (Theorem \ref{thm.Zar}), 
\item computations with {\sc Singular} (\cite{DGG}) for  some sporadic singularities in small characteristic.
\end{itemize}

Moreover, we give a characteristic-free proof of the results of Hefez and Hernandes \cite{HHA} in Section \ref{sec.hefher}. In addition to \cite{HHA} we give explicit lists of normal forms in Table \ref{tab.char0} and show that the classification in characteristic 0 holds also if the characteristic is bigger than the conductor of the semigroup. The difference between our method and that of \cite{HHA} is explained at the beginning of Section \ref{sec.hefher}.
\\

Our methodology is a combination of new theoretical results and symbolic computations, where the theoretical results 
lead to algorithms that allow explicit computations.  See in particular Proposition \ref{prop.modular}, where we  introduce the concepts of ``modular'' and ``full'' families,  which can be used, together with the semicontinuity of the modality,  to determine the modality in general.  
We refrain from formulating the algorithms abstractly but give detailed comments (see the proof of Lemma \ref{lem.a14},   the following lemmas,  and the remark after Lemma \ref{lem.a19}) which may present the general method even clearer. 
 Thus,  the concepts and methods of this paper may be a model for similar classification problems in local algebra and local algebraic geometry. \\

\begin{defn} \label{def.para}
Let $K$ be a field. A  {\em parameterization} of an (algebroid) plane curve singularity is a  (non-zero) morphism of local $K$-algebras
$$\psi:  K [[x,y]] \to K [[t]],  \ x\mapsto  x(t), \ \ y\mapsto y(t)$$ 
with $x(t), y(t) \in tK[[t]]$. For simplicity we write $\psi = (x(t),y(t))\in tK[[t]]^2$.
In concrete terms,
$$\big(x(t),y(t)\big)= \big( \sum_{i\geq1}a_it^i,  \sum_{i\geq1}b_it^i\big), \ a_i, b_i \in K.$$
The lowest non-vanishing term $\ord(\psi) := \min \{i\mid a_i \neq 0 \text{ or } b_i \neq 0\}$ is called the {\em order} or {\em multiplicity} of $\psi$.
 \end{defn}
 
The image of $\psi$ is the subring $K[[x(t),y(t)]]$ of $ K[[t]]$ and it can be shown that the kernel of $\psi$ is generated by a single, irreducible power series $f$ so that  $K[[x(t),y(t)]] \cong K[[x,y]]/\langle f\rangle$ is the local ring of an algebroid irreducible plane curve singularity, also called {\em plane branch}.  In the following we call $\psi$  also a {\em parametric plane curve (singularity)} or {\em parametric plane branch}.
We always assume that the parameterization  is {\em primitive}, i.e.  the {\em delta invariant}
satisfies
$$\delta(\psi):= \dim_{K} K[[t]]/K[[x(t),y(t)]]<\infty.$$

\begin{defn}\label{def.action}
Let $\mathcal{L}$ be the group of $K$-algebra
 automorphism of $K[[x,y]]$ and $\mathcal{R}$  the group of $K$-algebra
 automorphism of $K[[t]]$ and $\mathcal{A}=\mathcal{L} \rtimes \mathcal{R}$ the semi-direct product.  
 $\mathcal{A}$ is called the {\em left-right group} and acts on the space of parametric plane curve singularities
as 
$$ (l,r)\cdot \psi := l\circ \psi\circ r, \ l\in \mathcal{L},  \ r\in \mathcal{R}.$$
Two  parametric plane curves $\psi$ and $\phi$ are $\mathcal{A}-${\em equivalent}
(notation $\psi \sim_{\mathcal A} \phi$),  if they lie in the same orbit under the action of $\mathcal{A}$.
\end{defn}

Specifically, $\mathcal{A}$ acts on $tK[[t]]^2\setminus \{0\}$  as follows. We have
$$
\begin{array}{lll}
 l^{-1}(x) = ax +by + \sum_{|\alpha|>1}a_\alpha x^{\alpha_1}y^{\alpha_2}, \\
 l^{-1}(y) = cx +dy + \sum_{|\alpha|>1}b_\alpha x^{\alpha_1}y^{\alpha_2},
\end{array}
$$
 with $ad-bc\neq 0$ and $r(t) = et + \sum_{i>1}e_it^i$ with $e\neq0$. Then the  parameterization $\phi= r\circ \psi\circ l^{-1}$ is given by
$\phi(x)(t)= l^{-1}(x(r(t)))$ and $\phi(y)(t)=l^{-1}(y(r(t))).$

\begin{defn}\label{def.gamma}
Let $ \psi = (x(t), y(t))$  be a parameterized plane branch. The {\em semigroup} (of values) of $\psi$ is defined as 
$$\Gamma = \Gamma(\psi)=\{\ord_t(h)\mid h \in K[[x(t),y(t)]]\}.$$
 The {\em conductor} of $\Gamma$ or of $\psi$ is 
 $c(\Gamma): = c(\psi):= \min\{a \in \Gamma\mid a+b \in \Gamma \ \forall \ b \geq 0\}$.  It satisfies 
$c(\psi) = 2 \delta(\psi)$. The multiplicity $m(\Gamma)$ of $\Gamma$ is the minimal non-zero element in $\Gamma$.
\end{defn}

Our method of  classification makes essential use of the semigroup $\Gamma$ of $\psi$, which is an invariant of the $\mathcal A$-equivalence class\,\footnote{\,The semigroup encodes the equisingularity class of the plane branch.}.
We classify the small semigroups and fix the semigoup as the main invariant for the classification. We order the semigroups lexicographically and prove that if $\Gamma > \Gamma'$ and if the parameterization corresponding to $\Gamma'$ is not unimodal (resp. not simple) then the parameterization corresponding to $\Gamma$ 
is not unimodal (resp. simple). Using {\sc Singular} we find by explicit computations minimal semigroups corresponding to non-unimodal parameterizations, which excludes parameterizations with bigger semigroup.

The aim is to classify parameterizations of modality 1 (or unimodal parameterizations) with respect to $\mathcal A$-equivalence.
Informally, the modality of $\psi$ is the least number $m$ such that a small neighbourhood of $\psi$ in the space of primitive parameterizations  can be covered by a finite number of at most $m-$parameter families of orbits of the action of $\mathcal{A}$.
The precise definition is somewhat complicated and will be given in the next section  (Definition \ref{def.modal}).
\bigskip

 \section{Semicontinuity of the Modality}\label{sec.mod}
In order to define the modality for parameterizations and prove some properties, we will use methods and results from algebraic groups.  For this purpose, we assume from now on that the field $K$ is algebraically closed (this is assumed for our classification results anyway). We give a precise definition of  modality, prove its semicontinuity, and show that the modality is equal to the dimension of a full modular family. These are the main new results of this section.
\medskip

Since the infinite dimensional group $\mathcal{A}$ is not algebraic, the first step is to reduce the situation to a finite dimension by the following finite determinacy result. 
\begin{defn}[$\mA$-determinacy]\label{def.deter}
The parameterized plane curve singularity
$\psi \in tK[[t]]^2$ is called {\em $k-$determined w.r.t. $\mathcal A$} or {\em$ k - \mathcal A$-determined} if for every parameterization $\phi$ satisfying
$\psi \equiv \phi \text{ mod }  \langle t \rangle ^{k+1}K[[t]]^2$ we have  $\psi \sim_{\mathcal A} \phi$; $\psi$ is {\em finitely $\mathcal A$-determined} if it is $ k - \mathcal A$-determined for some positive integer $k$.
\end{defn}

The following determinacy bound and the semicontinuity of $\delta$ 
is proved in \cite{GP} for parameterizations of arbitrary (not necessary plane)  reduced curve singularities and for arbitrary fields.

\begin{prop} [\cite{GP}, Proposition 11]\label{prop.deter}
The parameterization $\psi$ is $(4\cdot\delta(\psi)-1)$-determined.
\end{prop}

It follows from Theorem \ref{thm.Zar} that for plane branches $\psi$ is even $2\delta(\psi)$-determined.\\

 While the exact bound is not so important, the finite determinacy of a parameterization is crucial
in order to define the modality  by using algebraic groups. 
We also need the semicontinuity of $\delta$ in a family of parameterizations.

\begin{defn}[Family of parameterizations]\label{def.fam}
Let $A=K[y_1,...,y_n]/I$ be an affine $K$-algebra of finite type and $Y  =V(I)\subset K^n$ the variety defined by the ideal $I$.\footnote{\,In this paper we consider only closed points, i.e., we identify $Y$ with $\Max(A)$, the set of maximal ideals in $A$.}  A {\em family of parameterizations} over $Y$ (or over $A$) is a morphism of $A$-algebras
$$\psi_A:  A [[x,y]] \to A [[t]],  \ x\mapsto  X(t), \ \ y\mapsto Y(t)$$ 
with $X(t), Y(t) \in tA[[t]]$, such that for fixed $y_0\in Y$ the induced map 
$$\psi_{A,y_0}: K [[x,y]] \to K [[t]],  (x,y) \mapsto (X(t) \text{ mod }\frak m_{y_0},Y(t)\text{ mod }\frak m_{y_0}) $$
with $\frak m_{y_0}\subset A$ the maximal ideal of $y_0\in Y$, is a parameterization of a plane curve singularity. We also call $\psi_A$ a {\em deformation} of $\psi_{A,y_0}$.

Specifically, let $(X(t),Y(t)) =  \big( \sum_{i\geq1}A_it^i,  \sum_{i\geq1}B_it^i\big)$ with  $A_i, B_i \in A$ such that  the images  $a_i, b_i$ of $A_i, B_i$ in $A/\frak m_{y_0}$ are not all zero. Then  $(x(t),y(t))=  \big( \sum_{i\geq1}a_it^i,  \sum_{i\geq1}b_it^i\big)$ is a parameterization in the sense of definition \ref{def.para}.
\end{defn}

\begin{thm}\cite[Theorem 18 (Semicontinuity of $\delta$)]{GP}\label{thm.semicont}
Let $\psi_A$ be a family of parameterizations over $Y$ and $y_0\in Y$. Then there exists an open neighbourhood $U$ of $y_0$ such that
$$\delta (\psi_{A,y}) \leq \delta (\psi_{A,y_0}) \text{ for all } y\in U.$$
In other words, for each $d>0$ the set 
$U_d = \{y\in Y\mid\delta(\psi_{A,y})\leq d\}$
is open in $Y$.
\end{thm}
The theorem is proved in \cite{GP} more generally for $A$ an arbitrary Noetherian ring and also for non-closed points in Spec($A$).\\

We have to consider $k$-jets of the involved power series.  For $l\in K[[t]]$ 
and $k>0$ let
$$j^k l:=\text{ image of } l \text{ in }J_k(t):=K[[{t}]]/\langle t \rangle^{k+1}$$
and for $r\in K[[x,y]]$ 
$$j^k r:=\text{ image of } r \text{ in }J_k(x,y):=K[[{x,y}]]/\langle x,y \rangle^{k+1}.$$
We identify elements of $J_k$ with their power series up to order $k$. 
\medskip

The $k$-jet of a parameterization
$\psi=(x(t),y(t))\in tK[[t]]^2$  \text{  is  } $$j^k\psi=(j^kx(t),j^ky(t))\in J_k(t)^2.$$
The $k$-jet $\mathcal{A}_k$  of $ \mathcal{A}$ is defined as follows.  An element $r$ in the right group $\mathcal{R}$ is uniquely determined by a pair of power series
$(r(x), r(y))\in \langle x,y\rangle K[[x,y]]^2$ with invertible linear part.
Then  the $k$-jet of $r$ is
$r_k:=(j^kr(x),j^k r(y))$ and $\mathcal{R}_k:=\{r_k\ |\ r\in  \mathcal{R}\}$ is the $k$-jet of $\mathcal{R}$. The definition of $\mathcal{L}_k :=\{l_k\ |\ l\in  \mathcal{L}\}$ is similar and we set 
$$\mathcal{A}_k:=\mathcal{L}_k \ltimes \mathcal{R}_k.$$ The algebraic group $\mathcal{A}_k$ is an affine $K$-variety, which acts algebraically on the affine $K$-variety $J_k(t)^2$ by
\begin{displaymath}
\begin{array}{ccccccc}
\mathcal{A}_k\times J_k(t)^2&\longrightarrow & J_k(t)^2, \quad
((r_k,\l_k),j^k \psi)&\mapsto &j^k(l_k\circ j^k\psi\circ r_k^{-1}).
\end{array}
\end{displaymath}
\medskip

We define now in general the modality of a point in an algebraic variety under the action of an algebraic group.
For any algebraic $K$-variety $X$ and algebraic $K$-group $G$ there exists a a Rosenlicht stratification $\{X_{i}, i=1,\ldots,s\}$.  $X$ is the finite disjoint union of the $G$-invariant locally closed algebraic subvarieties $X_{i}, i=1,\ldots,s$, such that $X_{i}/G$ is a geometric quotient with quotient morphism $p_{i}\colon X_{i}\to X_{i}/G$.  $X_1$ is open and dense in $X$. For $U\subset X$ open $(U\cap X_i)/G = p_i(U\cap X_i)$ is a constructible subset of
the $K$-variety $X_{i}/G$ with a well defined dimension, and $p_{i}\colon U\cap X_{i}\to X_{i}/G$ is a family of orbits of dimension $\dim (p_i(U\cap X_i))$.
We refer to \cite{GN} for details.

\begin{defn}\label{def.modal}\text{  (Modality) }
\begin{enumerate}
\item For $U$  an open subset of $X$ we define
$$G\text{-}\mathrm{mod}(U):=\max_{1\leq i\leq s}\{\dim \big (p_i(U\cap X_i)\big)\}.$$
For $x\in X$ we call 
$$G\text{-}\mathrm{mod}(x):=\min \{G\text{-}\mathrm{mod}(U)\ |\ U \text{ a neighbourhood of } x\}$$
the {\em $G$-modality} of $x$ (in $X$).
\item The {\em $\mathcal A$-modality} of a parameterization $\psi$, denoted by $\mathcal A\text{-mod}(\psi)$,  is defined as the 
$\mathcal A_k$-modality of $j^k\psi \in tJ_k(t)^2$ (in the sense of (1)) for sufficiently large $k$.  Hence
$$\mA\text{-mod}(\psi) = \mA_k\text{-mod}(U)=\max_i \{\dim (U\cap X_i)/\mA_k\}$$ 
for $U$ a sufficiently small neighbourhood of $j^k \psi$ and $\{X_i\}$ a Rosenlicht stratification of the jet-space $tJ_k(t)^2$. 

\item We call $\psi$ {\em simple, unimodal, bimodal} and {\em $r$-modal}  w.r.t. $\mathcal A$ if the $\mathcal{A}$-modality of  $\psi$ equals to 0, 1, 2 and $r$ respectively.
\item An integer $k$ is {\em sufficiently large} for $\psi$ w.r.t. $\mathcal A$ if there exists a neighbourhood $U$ of $j^k\psi$ in $tJ_k(t)^2$ s.t. every parameterization $\phi$ with $j^k\phi\in U$ is  $k$-determined w.r.t. $\mathcal A$.
\end{enumerate}
\end{defn}

It is shown in \cite[Corollary A.3]{GN} that the modality is independent of the Rosenlicht stratification. It is clear that $\mathcal A$-equivalent parameterizations have the same $\mathcal A$-modality.\\

Before we continue, we show that the jet spaces can be considered as families of parameterizations in the sense of Definition \ref{def.fam}.  This will be used in the sequel.

\begin{defn} [Jet space family]\label{def.jet}
(1) We identify the space of $l$-jets of parameterizations   
with an affine space  over $K$ by considering the coefficients of a parameterization as variables. Formally the identification is done by the map $tJ_l(t)^2 \to Y_l,$  
$$(\sum_{0<|\alpha|\leq l} a_{\alpha}{t}^{\alpha}, \sum_{0<|\beta|\leq l} b_{\beta}{t}^{\beta})\mapsto (x_{\alpha},y_\beta)_{0<|\alpha|,|\beta|\leq l},$$ 
where $\mathrm J_l:=K[x_\alpha,y_\beta]$ is the coordinate ring of $Y_l =\Max (\mathrm J_l)$.
Then
\begin{align*}
\begin{split}
\psi_{\mathrm J_l}:& \ \mathrm J_l[[x,y]] \to \mathrm J_l[[t]],\\
&\  x\mapsto X(t) = \sum_{0<|\alpha|\leq l} x_{\alpha}{t}^{\alpha},  \\
&\ y\mapsto Y(t)=\sum_{0<|\beta|\leq l} y_{\beta}{t}^{\beta},
\end{split}
\end{align*}
is a family of parameterizations over $Y_l$ in the sense of Definition \ref{def.fam}.
We call $\psi_{\mathrm J_l} = (\psi_{\mathrm J_{l,y}})_{y\in Y_l}$  the {\em jet space family (of parameterizations)}.
\smallskip

(2) Let $A$ be any affine $K$-algebra and $\psi_A:A[[x,y]] \to tA[[t]]$ be  a family of parameterizations over $Y=\Max(A)$ given by $(X(t),Y(t))$ (as in Definition \ref{def.fam}).  Then, for $y\in Y$,
$\psi_{A,y}$ is the image of $(X(t),Y(t))$ in 
$(tA[[t]]/\frak m_ytA[[t]])^2=tK[[t]]^2$. Taking $l$-jets of the $\psi_{A,y}$ we get a morphism 
$$\Psi_l:Y\longrightarrow tJ_l(t)^2 = Y_l, \ y\mapsto j^l\psi_{A,y},$$
the {\em l-jet map} of $\psi_A$.
\end{defn}

\begin{prop}\label{prop.suff}
Let $\psi$ be a parameterization of a plane curve singularity.
Then every $k\geq 4\cdot\delta(\psi)$ is sufficiently large for $\psi$ w.r.t. $\mathcal A$. Moreover,  $\mathcal A\text{-mod}(\psi)$ is independent of the sufficiently large $k$.
\end{prop}

\begin{proof}
Consider the jet space family $\psi_{\mathrm J_l}:\mathrm J_l[[x,y]] \to  \mathrm J_l[[t]]$. By the upper semi-continuity of $\delta$ (Theorem \ref{thm.semicont}),  for each $l,d$ the subset
$U_{l,d}:=\{y\in Y_l\mid\delta(\psi_{\mathrm J_l,y})\leq d\}$
is open in $Y_l$.  For each parameterization $\phi$
the $l$-jet $j^l\phi$ corresponds to a point $y_l\in Y_l$ such that  
$j^l\phi = \psi_{\mathrm J_l,y_l}$. In particular, the set
$$U_l:=\{j^l \phi\in tJ_l(t)^2\mid\delta(j^l\phi)\leq \delta (j^l \psi)\}$$
is an open neighbourhood of $j^l(\psi)$. 

Now let $k\geq 4\delta(\psi)$. Then each parameterization $\phi$ with $\delta(\phi)\leq\delta(\psi)$ is  $k$-determined w.r.t. $\mathcal A$ by Proposition \ref{prop.deter} and satisfies $\delta(\phi) = \delta(j^k\phi)$.  It follows that
 every parameterization $\phi$ with $j^k\phi\in U_k$ is  $k$-determined w.r.t. $\mathcal A$.
 This means that $k$ is sufficiently large for $\psi$ w.r.t. $\mathcal A$. 
 
 The last statement can be proved in a similar way as in the proof of \cite[Proposition 2.6]{GN} (for right and contact equivalence of hypersurface singularities). We omit the details.
\end{proof}

We can now prove the main result of this section, the semicontinuity of the $\mathcal A$-modality in a family of parameterizations.

\begin{thm}[Semicontinuity of $\mA$-modality]\label{thm.semicontmod}
Let $A$ be an affine $K$-algebra and $\psi_A$ a family of parameterizations over $Y=\Max(A)$ (as in Definition \ref{def.fam}) and let $y_0\in Y$. Then there exists an open neighbourhood $U\subset Y$ of $y_0$ such that
$$\mathcal A\text{\em -mod}(\psi_{A,y}) \leq \mathcal A\text{\em -mod} (\psi_{A,y_0}) \text{ for all } y\in U.$$
It follows that for each $m\geq 0$  the set
$$U_{m}:=\{y\in Y\mid\mathcal A\text{\em -mod}(\psi_{\mathcal A,y})\leq m\}$$
is open in $Y$. 
\end{thm}

\begin{proof}
(1) Consider first for $l\geq 1$  the jet space family $\psi_{\mathrm J_l}$  over $Y_l$.  We show that for each $m\geq 0$ the subset
$U_{l,m}:=\{y\in Y_l\mid\mathcal A\text{\em -mod}(\psi_{\mathrm J_l,y})\leq m\}$
is open in $Y_l$. 

To see this,
take an open neighbourhood $V$ of $y=j^l\psi=\psi_{\mathrm J_l,y}$ in $Y_l = tJ_l(t)^2$ such that $\mathcal A_l\text{-}\mathrm{mod}(V)=\mathcal A_l\text{-}\mathrm{mod}(j^l\psi)$ (cf. Definition \ref{def.modal}). Let
$$\tilde U:=\{j^l\phi\in tJ_l(t)^2  \mid \delta(j^l\phi)\leq \delta(j^l\psi)\}.$$
By Theorem \ref{thm.semicont} $\tilde U$ is open and hence $U:= V\cap \tilde U$ is an open neighbourhood of $j^l \psi$ in $tJ_l(t)^2$. 
 For $j^l \phi \in U$ we have 
$$\mathcal A_l\text{-}\mathrm{mod}(j^l\phi)\leq \mathcal A_l\text{-}\mathrm{mod}(V)=\mathcal A_l\text{-}\mathrm{mod}(j^l\psi).$$
Hence $U\subset U_{l,m}$ and $U_{l,m}$ is open in $tJ_l(t)^2$.

(2) For the general case let $\Psi_k:Y\longrightarrow tJ_k(t)^2 = Y_k, \ y\mapsto j^k\psi_{A,y}$ be the $k$-jet map of $\psi_A$ (Definition \ref{def.jet}).
By (1) the set $U_k= \{j^k \phi \mid \mathcal A_k\text{-mod}(j^k\phi) \leq \mathcal A_k\text{-mod} (j^k\psi_{A,y_0}\}$ is open in $Y_k$.
If $k\geq 4\delta(\psi_{A,y_0})$, then $k$ is sufficiently  large  for $\psi_{A,y_0}$ by Proposition \ref{prop.suff} and there exist an open neigbourhood $V \subset Y_k$ of $j^k\psi_{A,y_0}$ such that $\psi_{A,y} \sim_\mathcal A j^k\psi_{A,y}$ for $y\in V$.
Then $U:=\Psi_k^{-1}(U_k\cap V)$ is an open neighbourhood of $y_0$ in $Y$ and we get for $y\in U$
$$  \mathcal A\text{-mod}(\psi_{A,y}) = \mathcal A_k\text{-mod}(j^k\psi_{A,y}) \leq
\mathcal A_k\text{-mod}(j^k\psi_{A,y_0}) = \mathcal A\text{-mod}(\psi_{A,y_0}),$$
proving the theorem.
\end{proof}

The following interpretation of the modality as the dimension of some "full modular family" is quite useful and used in our classification. Informally speaking, a family of parameterizations is "modular" if it does not contain a positive dimensional trivial subfamily and it is "full" at a given parameterization $\psi$ if it contains all deformations of $\psi$  up to $\mA$-equivalence.

\begin{defn}[Modular and full family]\label{def.modular}
Let  $\psi_A$  be a family of parameterizations over $Y=\Max(A)$ (as in Definition \ref{def.fam}).

(1)  $\psi_A$ is called {\em ($\mathcal A$-)modular} if for each $y\in Y$ there are only finitely many $y'\in Y$ such that $\psi_{A,y}$ is $\mathcal A$-equivalent to $\psi_{A,{y'}}$. 

(2)  $\psi_A$ is called {\em ($\mathcal A$-)full at $y_0 \in Y$}, or  {\em full at $\psi =\psi_{A,{y_0}}$}, if  the following holds:
Let $k$ be sufficiently large for $\psi$
 and let $U$ be an open neighbourhood of $j^k\psi$ in $tJ_k(t)^2$ such that $\mathcal A$-{mod}$(j^k\psi)= \mathcal A$-{mod}$(U)$. Then  for each  $j^k \phi \in U$  there exists a $y\in Y$ with $j^k\phi  \sim_\mathcal A \psi_{A,{y}}$. 
 
(3) We say that  $\psi_A$ is {\em $\Gamma$-full at $y_0$ }if  for each  $j^k \phi \in U$ s.t.  $\Gamma (\psi) = \Gamma (j^k\phi)$ there exists a $y\in Y$ with $j^k\phi  \sim_\mathcal A \psi_{A,{y}}$.  
\end{defn}

\begin{prop}\label{prop.modular}
Let $\psi$ be a parameterized plane branch and $\psi_A$  a family of parameterizations over $Y$  with $Y$ irreducible\,\footnote{\,The proof shows that  for (2) we do not need that $Y$ is irreducible and for (1) we need only that $\dim \Psi_k(Y) = \dim \Psi_k(Y)\cap U$, which holds if $U$ is dense (e.g., if $Y$ is irreducible).
} and $y_0 \in Y$ such that $\psi \sim_{\mathcal A} \psi_{A,{y_0}}$. 
\begin{enumerate}
\item If $\psi_\mA$ is modular, then 
$\mathcal A\text{\em -mod}(\psi) \geq \dim Y.$
\item If $\psi_\mA$ is full at $y_0$, then
$\mathcal A\text{\em -mod}(\psi) \leq \dim Y.$
\item If $\psi_\mA$ is modular and full at $y_0$ then
$\mathcal A\text{\em -mod}(\psi)= \dim Y.$
\item Let $\psi_\mA$ be modular and $\Gamma$-full at $y_0$. Suppose that 
$\Gamma(\psi)$ has $\leq 3$ generators and
for any parameterization $\psi'$ with $\Gamma(\psi') < \Gamma(\psi)$ (lexicographically) holds
$\mathcal A\text{\em -mod}(\psi') \leq \dim Y$. Then 
$\mathcal A\text{\em -mod}(\psi)= \dim Y.$
\end{enumerate}
\end{prop}

\begin{proof}
Consider the jet space family $\psi_{\mathrm J_k}:\mathrm J_k[[x,y]] \to  \mathrm J_k[[t]]$ and the $k-$jet morphism of $\psi_A$,
$\Psi_k:Y\to tJ_k(t)^2, \ y\mapsto j^k\psi_{A,y}$ (Definition \ref{def.jet}).
Let $k$ be sufficiently large for  $\psi$ and 
$U$ an open neighbourhood of $j^k\psi$ in $tJ_k(t)^2$ such that $\mathcal A\text{-mod}(\psi)=\mathcal A$-mod$(j^k(\psi))= \mathcal A$-mod$(U)$.
Since $Y$ is irreducible, $U$ is dense in $Y$. 

(1) Since $\psi_A$ is modular, the map $\Psi_k$ has finite fibers, because for any two points $y, y'$ in a fiber  the parameterizations $\psi_{A,y}$ and $\psi_{A,y'}$ are isomorphic (they have the same $k$-jet). It follows that $\dim Y = \dim \Psi_k(Y)$ and  there are only finitely many parameterizations in $ \Psi_k(Y) \subset  tJ_k(t)^2$ that are $\mA$-equivalent to a given one. Let $X_1$ be the open and dense stratum of a Rosenlicht stratification of $ tJ_k(t)^2$. Then $\mathcal A$-mod$(j^k \psi_{A,{y_0}})= \mathcal A$-mod$(U) \geq \dim U\cap X_1/\mA_k$. 
The quotient map $p_1: U\cap X_1 \to U\cap X_1/\mA_k$ restricted to $\Psi_k(Y)\cap U \cap X_1$ has finite fibers and since $X_1$ and $U$ are dense, $\dim Y = \dim \Psi_k(Y) = \dim \Psi_k(Y)\cap U \cap X_1 =  \dim \Psi_k(Y)\cap U \cap X_1/\mA_k \leq \dim U\cap X_1/\mA_k \leq \mA$-mod$(\psi)$.

(2) Consider again  $\Psi_k: Y \to tJ_k(t)^2$ and let $X_i$ be a stratum of the Rosenlicht stratification such that $\mathcal A$-mod$(U) = \dim U\cap X_i/\mA_k$. Let $p_i: U\cap X_i \to U\cap X_i/\mA_k $ be the quotient map.
Since $\psi_A$ is full, for each class $[j^k\phi ]\in U\cap X_i/\mA_k$ the fiber 
$(p_i \circ \Psi_k)^{-1}([j^k\phi])$ is non-empty. Hence $\dim Y \geq \dim U\cap X_i/\mA_k = \mathcal A$-mod$(\psi)$.

(3) follows from (1) and (2).

(4)  If $\mathcal A\text{-mod}(\psi) > \dim Y$,  $\psi$ could be deformed to a
$\psi'$  with $\mathcal A\text{-mod}(\psi') = \mathcal A\text{-mod}(\psi)$
and with
$\Gamma(\psi') \neq \Gamma(\psi)$ since $\psi_\mA$ is $\Gamma$-full. By the following proposition $\Gamma(\psi') <\Gamma(\psi)$ and thus 
$\mathcal A\text{-mod}(\psi') \leq \dim Y$, a contradiction.
\end{proof}

\begin{prop}[Semicontinuity of $\Gamma$]\cite[Lemma 14]{MP1}\label{prop.semicontG}
Let $\psi(t)$ be a parameterization with semigroup $\Gamma$ generated by
at most $3$ elements and $\psi_A$ a deformation of $\psi = \psi_{A,y_0}$. Then the semigroups  of
$\psi_{A,y}(t)$ for $y$ in a neighbourhood of $y_0$ are lexicographically smaller or equal to $\Gamma$.
\end{prop}
\begin{rem} \label{rem.sem}
The semicontinuity of semigroups with more than three generators is, to the best of our knowledge,  an open problem.
\end{rem}
\bigskip


\section{Methods of proof and overview of results}
\label{sec.meth}
For the whole section let $K$ be an algebraically closed field
 of characteristic $p\geq 0$. 
An important ingredient of the proof is a generalization of Zariski's ``short parameterization'' to positive characteristic.
Let $(x(t), y(t))$ be a parameterized plane branch  defined over $K$ with semigroup $\Gamma$ and conductor $c$.  Let 
$$a:= \ord_t(x(t)) < b:= \ord_t(y(t)).$$

In \cite[Chapter III, Proposition 1.2]{ZO} Zariski proved for complex analytic branches  the existence  of a ``short parameterization''  of the form $x(t)=t^a$,  $ y(t) = t^b+\sum_{i \notin \Gamma}b_it^i$ (up to $\mathcal A$-equivalence). For algebroid
 branches over an algebraically closed field of any characteristic 
we extend (and improve) Zariski's proof to get ``short parameterizations''of different  kinds (depending on $p$) that will be useful in our classification (by definition, $p=0$ implies that $p$ does not divide any integer).

\begin{thm}[Short parameterization in arbitrary characteristic] \label{thm.Zar} Let $(x(t), y(t))$ be a parameterization  over $K$  with semigroup $\Gamma$ and char$(K)=p\geq 0$.
\begin{enumerate}
\item In general $(x(t), y(t))\sim_{\mathcal A}(t^a+\sum_{i \notin \Gamma}a_it^i,t^b+\sum_{i \notin \Gamma}b_it^i)$  and $a\nmid b$. \\
If  $p \nmid a $ then $(x(t), y(t)) \sim_{\mathcal A}(t^a, t^b+\sum_{i \notin \Gamma}b_it^i)$.\\
If  $p \nmid b $ then $(x(t), y(t)) \sim_{\mathcal A}(t^a+\sum_{i \notin \Gamma}a_it^i, t^b)$.
\item Let $a>2$.  If $p \nmid b$ then 
$(x(t), y(t)) \sim _{\mathcal{A}}(t^a+\sum_{i \notin \Gamma}a_it^i,t^b+\sum_{i \notin \Gamma, i<c-1}b_it^i).$
Similarly, if $p \nmid a$ then $(x(t), y(t)) \sim _{\mathcal{A}} (t^a+\sum_{i \notin \Gamma, i<c-1}a_it^i, t^b+\sum_{i \notin \Gamma}b_it^i).$
\item If $a>2$, $p\nmid a$ and $p\nmid b$ then\,\footnote{\,It is important for the applications that the term of order ${c-1}$ can be eliminated (this observation seems to be new). }
$$\begin{array}{clll}
(x(t), y(t)) &\sim _{\mathcal{A}}  (t^a, t^b+\sum_{i \notin \Gamma,i<c-1}b_it^i)\\
& \sim _{\mathcal{A}} (t^a+\sum_{i \notin \Gamma, i<c-1}a_it^i, t^b).\end{array}$$
\end{enumerate}
\end{thm}

\begin{rem} \label{rm.Zar}
(1) The families on the r.h.s. with parameters $a_i, b_i$ are $\Gamma$-full in the sense of Definition \ref{def.modular} (but not necessarily modular).

(2) For $a\leq 2$ we have a more precise result. 
$a=1$ means that the parameterized branch is nonsingular with $(x(t),y(t)) \sim _{\mathcal{A}} (t,0)$.  If $a=2$ then $b=2k+1$ is odd 
and if $p\neq 2$ then $(x(t),y(t)) \sim _{\mathcal{A}} (t^2,t^{2k+1})$. If 
$p= 2$ then $(x(t),y(t))$ is $\mathcal{A}$-equivalent to either 
$(t^2,t^{2k+1})$ or to $(t^2+t^{2m+1},t^{2k+1})$, $0<m<k$  (cf. \cite[Theorem 3]{MP1}).  In particular,  the parameterization is 0-modal (simple) if $a\leq 2$. 

(3) Since the r.h.s of the theorem involves only exponents of $t$, which are not in $\Gamma$, it follows that the $\mathcal{A}$-determinacy of a plane branch is bounded by the conductor $c(\Gamma)$.
\end{rem}

\begin{proof}
(1) We note that \cite[Lemma 1.1]{ZO} holds in any characteristic. The transformations of Zariski \cite[Proof of Proposition 1.2]{ZO} can be applied to $x(t)$ and $y(t)$ in any characteristic, showing  that $(x(t),y(t))$ is in general $\mathcal A$-equivalent to the r.h.s.  of (1).
E.g., the transformation to eliminate $a_\nu t^\nu$  from $x(t) = t^a+...+a_\nu t^\nu + ...$, with $a<\nu<c$ minimal s.t. $\nu \in \Gamma$, is of the form $y\mapsto y$ and $x\mapsto x-h $ with $h\in K[[x(t),y(t)]]$, $\ord_t(h)=\nu$. This transformation does not change the terms $a_it^i$ in $x(t)$ with $i<\nu$. Similar for $y(t)$.

If $p\nmid a$ then the $a$-th  root  $u(t)$ of  $1+ \sum_{i \notin \Gamma}a_it^{i-a}$ exists and  $t \mapsto tu(t)$ is an automorphism of $K[[t]]$.  The inverse transformation $\psi$ maps $x(t)$ to $t^a$ and $\psi(y(t))$ 
can be transformed to the required form by Zariski's argument from the previous step.  An analogous argument applies if $p\nmid b$.

(2)
Let $p \nmid b$. Since $a>2$ we have $b<c-1$.
If $b_{c-1}\neq 0$ we consider the transformation $\varphi(t)=t-\frac{1}{b}b_{c-1} t^{c-b}$,  $c-b\geq 2$. Checking exponents we obtain 
$$\begin{array}{clll}
\varphi(x(t))=&(t-\frac{1}{b}b_{c-1} t^{c-b})^a+\ldots=t^a-\frac{a}{b}b_{c-1}t^{c-1-(b-a)}+\ldots, \\
\varphi(y(t))=&t^b+\sum_{i \notin \Gamma, i<c-1}b_it^i + \text{ terms of order} \geq c.
\end{array}$$
 Now apply Zariski's transformations from (1)  to get the desired form.

If  $p \nmid a$ and $a_{c-1}\neq 0$ we use the transformation $t\mapsto t-\frac{1}{a}a_{c-1} t^{c-a}$  to get the result.

(3) By (1) we can assume that $(x(t), y(t)) = (t^a, t^b+\sum_{i \notin \Gamma}b_it^i)$.
$\Gamma$ is symmetric and $b-a\notin \Gamma$. Therefore $c-1-(b-a)\in \Gamma$. Using $\varphi$ from (2) we have  $\varphi(t^a)=t^a-\frac{a}{b}b_{c-1}t^{c-1-(b-a)}+$ terms of order $ > c-1-(b-a)$ and we get as in (2)
 \begin{center}
 $(\varphi(t^a),\varphi(y(t))) \sim _{\mathcal{A}} (t^a+\sum_{i \notin \Gamma, i>c-1-b+a}a_it^i,t^b+\sum_{i \notin \Gamma, i<c-1}b_it^i )$.
 \end{center}
 Now  define the map $\psi$ by 
 \begin{center}
 $\psi(t^a+\sum_{i \notin \Gamma, i>c-1-b+a}a_it^i)=t^a$. 
 \end{center}
 Since $\psi(t)=t+\text{ terms of order } > c-1-(b-a)$ we obtain 
 \begin{center}
 $\psi(t^b+\sum_{i \notin \Gamma, i<c-1}b_it^i )=t^b+\sum_{i \notin \Gamma, i<c-1}b_it^i + \text{ terms of order} \geq c$,
 \end{center}
where the terms of order $\geq c$ can be eliminated by (1). This proves the first claim.
The second claim follows similarly.
\end{proof}
\bigskip

Let us now give an overview of some methods and the results of this paper.
The proofs will be given in the subsequent sections.

\begin{itemize}
\item As the main invariant we fix the semigroup  
$$\Gamma =\Gamma(\psi) =\langle n_1,n_2,...,n_r\rangle, 1< n_1 < n_2, ...<n_r,$$
 given by a minimal system of generators $\{n_i\}$, with $n_1$ the multiplicity.  Some semigroups depend on a parameter
$k\in \N$, such as $\langle 4,10,2k+11\rangle$,  the others, like  $\langle 4,9\rangle$, are called ``sporadic''.
\medskip

We order the semigroups lexicographically. By Remark \ref{rm.Zar} we know that 
 in case of $n_1\leq 2$  the parameterization is simple.
In section 6 we show that the parameterization is not unimodal if the multiplicity is $>5$ in characteristic $\neq 2$,
cf.  Proposition \ref{notuni}.
\item For the following semigroups $\Gamma$ we check the unimodality by explicit construction of 1-dimensional $\Gamma$-full modular families and apply Proposition \ref{prop.modular}.  

We set $k_0 := \min \{k \text{ such that } p\mid 2k+7\}$, i.e.,  
\begin{align} \tag{*}\label{no.1}
\begin{split}
       k_0=&\frac{p-7}{2} \text{ for } p\geq 13,\\
       k_0=&13 \text{ for } p=11,                \\      
       k_0=&7  \text{ for } p=7,     \text{ and }       \\         
        k_0=&4 \text{ for }  p=5 \text{ and } p=3. 
\end{split}
\end{align}
Consider the following {\bf  list of semigroups of multiplicity 3,4,5:}
\begin{itemize}
\item $\Gamma = \langle 3,k\rangle$, \ $k \geq p+9$,  $3\nmid k$,
\item $\Gamma = \langle 4,6,2k+7\rangle$, \ $k\geq k_0$,
\item $\Gamma = \langle 4,9\rangle$, 
\item $\Gamma = \langle 4,10,2k+11\rangle$, \ $k\geq 5$,
\item $\Gamma = \langle 4,11\rangle$, 
\item $\Gamma = \langle 5,6\rangle$, 
\item $\Gamma = \langle 5,7\rangle$, 
\item $\Gamma = \langle 5,8\rangle$. 
\end{itemize}
\medskip

In \cite{MP1} it was shown:
\item parameterizations of multiplicity 3.\\
Let $\Gamma=\langle 3,k\rangle$ and\,\footnote{\,The condition on $k$ is equivalent to the cases $k=3\bar k+1$, $l=3 \bar l +2$ and  
$k=3\bar k+2$, $l=3 \bar l +1$, which was used in \cite{MP1}.} $(t^3,t^k+t^l+\sum_{i>l}a_it^i)$, $k<l$, $kl\equiv 2 \text{ mod } 3$, $l\leq 2k-9$, be a parameterization then
\begin{enumerate}
\item If $p\nmid l-k$ then  $(t^3,t^k+t^l+\sum_{i>l}a_it^i)$ is equivalent to $(t^3,t^k+t^l)$, (lemma 4 and 7).
\item If $p\mid l-k$ then  $(t^3,t^k+t^l+\sum_{i>l}a_it^i)$ is not simple (lemma 5 and 8).\\
This implies the following (using deformation arguments\,\footnote{\,Deform $(t^3,t^k+t^l+\sum_{i>l}a_it^i)$ to
$(t^3,t^k+\alpha t^{k+p}+t^l+\sum_{i>l}a_it^i)$ in case of (3) respectively to $(t^3,\alpha t^{p+9}+\beta t^{2p+9}+t^k+t^l+\sum_{i>l}a_it^i)$ in case of (4).}
and semicontinuity of modality):
\item If $l\geq k+p$ (and therefore $k\geq p+9$) and $k(k+p)\equiv 2 \text{ mod } 3$ then $(t^3,t^k+t^l+\sum_{i>l}a_it^i)$ is not 
simple\footnote{\,Note that $l \geq k+2p$ implies $k\geq 2p+9$.}. 
\item The parameterization $(t^3,t^{p+9}+t^{2p+9})$ is not simple\footnote{\,Since $p^2\equiv 1 \text{ mod } 3$ this implies $(p+9)(2p+9)\equiv 2 \text{ mod } 3$.}. This implies $(t^3,t^k+t^l+\sum_{i>l}a_it^i)$ is not simple if $k\geq 2p+9$.
\end{enumerate}
We obtain that parameterizations with semigroup $\langle 3,k\rangle$ are not simple if $k\geq 2p+9$.
If $k<2p+9$ then the parameterization $(t^3,t^k+t^l+\sum_{i>l}a_it^i)$ is not simple iff $k\geq p+9$, $l\geq k+p$ and
$k(k+p)\equiv 2 \text{ mod } 3$.

Note that (3) and (4) replace Corollary 5 from \cite{MP1}, which is not fully correct.
\item 
parameterizations of multiplicity 4 with semigroup $\Gamma < \langle 4,6,2k_0+7\rangle$ (lexicographically) or $\Gamma=\langle 4,7\rangle$  are simple.
\medskip

We shall prove:
\item parameterizations of multiplicity 4 with semigroup $\Gamma >\langle 4,11\rangle$ are not unimodal (Proposition \ref{notuni}).
\item parameterizations with semigroup $\Gamma >\langle 5,8\rangle$ are not unimodal (Proposition \ref{notuni}).
\item The classification for $p=2,3$ is done in Section \ref{sec.class3}.
\item The results are presented in tables of normal forms of parameterizations.
\end{itemize}

\begin{exmp}
The following (not expected) example shows that the modality of a parameterization with semigroup $\langle 3,k\rangle$ may increase while $k$ decreases, if the characteristic is positive.
Let $p=7$ and consider $(x(t),y(t))=(t^3,t^k+t^{k+p}+\sum_{i>k+p}a_it^i)$
with semigroup  $\langle 3,k\rangle$  and conductor $2k-2$.

If $k=17$ then $3\mid k+p =24$ and  $(x(t),y(t))\sim_\mA (t^3,t^k+\sum_{i>k+p}a_it^i)$. Let $i\leq k+2p$ be minimal with $a_i\neq 0$. If $i < k+2p$ then 
$(x(t),y(t))\sim_\mA (t^3,t^k+t^i)$ by (1) above, if $i \geq k+2p$ then $(x(t),y(t))\sim_\mA (t^3,t^k)$ by Theorem \ref{thm.Zar}. In any case, $(x(t),y(t))$ is simple.

If $k=16$ then $k+p=23$ with $k(k+p)\equiv 2 \text{ mod } 3$ and $(x(t),y(t))$ is not simple by (3) above. Indeed it is unimodal with normal form 
$(t^3,t^{16}+t^{23} + at^{26})$.
\end{exmp}
From the above remarks and the semicontinuity of $\Gamma$ (Proposition \ref{prop.semicontG}) we get the following Lemma.

\begin{lem}
The above  semigroups of multiplicity 3,4,5 are the only candidates for semigroups of unimodal parameterizations.
\end{lem}
\bigskip

The proof of the following theorem will be given in the subsequent sections.

\begin{thm} \label{thm.main}
Let $(x(t), y(t))$ be a parameterized plane branch 
defined over an algebraically closed field $K$ of characteristic $p>0$. 
The parameterization is $\mathcal A$-unimodal if and only if it is contained in one of the following tables.
\begin{enumerate}
\item If $p\geq 5$  in Table \ref{tab.p>4,m<5} and Table \ref{tab.p>4,m=5}.
\item If $p=3$  in Table \ref{tab.p=3}.
\item If $p=2$  in Table \ref{tab.p=2}. 
\end{enumerate}
\end{thm}

  \begin{table} [h!]   
   \caption{Unimodal Parametric Plane Curve Singularities, $p\geq 5$, 
   multiplicity $\leq 4$ (3 series, 13 sporadic singularities, $a\in K$)}\label{tab.p>4,m<5}
       \begin{center}
   \begin{tabular}{| l | l |  c | lc | }
      \hline
       Normal form              & Semigroup   
       \\     \hline   
      \raisebox{-2 pt}{$(t^3,t^k+t^l+at^{l+3})$ with }
      &  \raisebox{-2 pt}{$\langle 3,k \rangle $ with $k\geq 2p+9$}\\
      $k\leq l \leq 2k-9,$ $k \cdot l \equiv 2 \text{ mod }3$ 
      &  or \big ($p+9\leq k <2p+9$   and\\
      \raisebox{1 pt}{$(t^3,t^{k})$ }  &  \raisebox{1 pt}{$l\geq k+p$, $k(k+p)\equiv 2 \text{ mod } 3$\big )}
           \\  \hline 
      \raisebox{-2 pt}{$(t^4,t^6+t^{2k+1}+at^{2k+3})$ } 
      & \raisebox{-2 pt}{$\langle 4,6,2k+7 \rangle$ with
        $k\geq k_0$   (\ref{no.1})}        
       \raisebox{-6 pt}{\text{}} \\
       \hline  
       
       \raisebox{-2 pt}{$(t^4, t^9+t^{10}+at^{11})$, $a\neq \frac{19}{18}$}      &  \\
       \raisebox{-1 pt}{$(t^4, t^9+t^{10}+\frac{19}{18}t^{11}+at^{15})$ }                     &  \\
       $(t^4, t^9+t^{11})$                     & $\langle 4,9 \rangle$ \\
       $(t^4, t^9+t^{15})$                     &  \\
       $(t^4, t^9+t^{19})$                 &  \\
       \raisebox{1 pt}{$(t^4, t^9)$ }                &  \\\hline

        \raisebox{-2 pt}{$p\neq 7$ }
       &  \raisebox{-2 pt}{$\langle 4,10,2k+11 \rangle$,  $k<\frac{p-11}{2}$ if $p\geq 23$,} \\
       $(t^4, t^{10}+t^{2k+1}+at^{2k+3})$                   
       & $k<\frac{p+9}{2}$ if $5\leq p\leq 19$ 
        \raisebox{-6 pt}{\text{}} \\
                 \hline

       \raisebox{-2 pt}{$p \neq 11$ }                                          &\\     
       $(t^4, t^{11}+t^{13}+at^{14})$               &  \\
        $(t^4, t^{11}+t^{14}+at^{17})$, $a\neq \frac{25}{22}$               &  \\
        $(t^4, t^{11}+t^{14}+\frac{25}{22}t^{17}+at^{21})$   
        &  $\langle 4,11 \rangle$                   \\
       $(t^4, t^{11}+t^{17})$                     &   \\
       $(t^4, t^{11}+t^{21})$                     &  \\
       $(t^4, t^{11}+t^{25})$                     &  \\
      \raisebox{1 pt}{$(t^4, t^{11})$ }                               &  \\ \hline
    \end{tabular}
       \end{center}
    \end{table}
\newpage    
     \begin{table} [h!]  
  \caption{Unimodal Parametric Plane Curve Singularities,  $p\geq 5$, multiplicity $5$ (15 sporadic singularities, $a\in K$).}  
  \label{tab.p>4,m=5}
        \begin{center}
    \begin{tabular}{| l | l |c|c|}
      \hline
       Normal form              & Semigroup    \\\hline
        $p\neq 5$                                 &\\
         $(t^5, t^{6}+ t^{8}+at^{9})$                &  \\
        $(t^5, t^6+t^9)$                    &$\langle 5,6 \rangle$ \\
        $(t^5, t^{6}+ t^{14})$  &  \\
        $(t^5, t^{6})$   &  \\\hline

       $p\neq 5,7$                                     &\\
       $(t^5, t^{7}+t^{8}+at^{11})$                &  \\
        $(t^5, t^{7}+ t^{11}+at^{13})$                    &$\langle 5,7 \rangle$ \\
        $(t^5, t^{7}+ t^{13})$  &  \\
        $(t^5, t^{7}+ t^{18})$  &  \\
        $(t^5, t^{7})$   &  \\\hline

        $p\neq 5$                                     &\\
        $(t^5, t^{8}+t^{9}+at^{12})$                &  \\
        $(t^5, t^{8}+ t^{12}+at^{14})$                    &$\langle 5,8 \rangle$ \\
        $(t^5, t^{8}+ t^{14}+at^{17})$  &  \\
        $(t^5, t^{8}+ t^{17})$  &  \\
        $(t^5, t^{8}+ t^{22})$  &  \\
        $(t^5, t^{8})$   &  \\ \hline

    \end{tabular}
   \end{center}
    \end{table}
\bigskip

\section{Classification in characteristic $\geq 5$}
\label{sec.class5}
Parameterizations with multiplicity 2 are simple (Remark \ref{rm.Zar}) and parameterizations with multiplicity $\geq 6$ are at least bimodal if   $p=$ char$(K) \neq 2$ (Proposition \ref{notuni}). We thus have to consider only the multiplicities 3,4,5 if   $p\neq 2$. In this section we assume $p\geq 5$, except we say otherwise.\\

We give a rough description of our procedure (in every characteristic). We start from a short  parameterization (Theorem \ref{thm.Zar}) $\psi_{\bf a}(t)= \psi(t) =(x(t),y(t))$ with a given semigroup that depends on certain parameters ${\bf a}=(a_i)$. Then we construct explicit automorphisms
$r: t\mapsto \varphi(t)$ in the right group $\mathcal{R} =Aut_KK[[t]]$ and automorphisms 
$l :(x,y)\mapsto (H(x,y), L(x,y))$  in the left group  $\mathcal{L} = Aut_K K[[x,y]]$ such that 
\begin{align}\tag{**}\label{no.2}
\psi'(t)&:=((H,L)\circ \psi\circ \varphi) (t) =(H(x(\varphi(t)),y(\varphi(t)),L(x(\varphi(t)),y(\varphi(t)))
\end{align}
is $\mathcal A$-equivalent to  $\psi(t)$ and has less parameters.  In the following computations we often replace $\varphi$ by $\varphi^{-1}$ and write (\ref{no.2}) as an "Ansatz" 
\begin{align*}
&x'(\varphi(t)) = H_x(x(t),y(t)), \\
&y'(\varphi(t)) = H_y(x(t),y(t)),
\end{align*}
with $\varphi$, $H_x$ and $H_y$ having unknowns as coefficients  that have to be determined.

If we reach at a parameterization $\psi_{a}(t)$ depending  only on one parameter $a$ such that  $\psi_a(t) \sim_{\mathcal A}\psi_{a'}(t)$  for only finitely many pairs $(a,a')$, the family is modular and $\Gamma$-full. Using the semicontinuity of the semigroup under deformations, we conclude that the family is even full and then the parameterization  is $\mathcal A$-unimodal by Proposition \ref{prop.modular}.

If, on the other hand, the number of parameters cannot be reduced to 1, the parameterization is not unimodal (if all parameters can be eliminated, the parameterization is 0-modal, i.e., simple).
For some  computations we use {\sc Singular} \cite{DGG}.

 \subsection{Multiplicity 3} $\text{}$\\

 We assume that the characteristic $p$ of $K$ is different from $3$ ($p=2$ is allowed).\\
 We have the semigroups  $\langle 3, k\rangle$,  $k\geq 4$ and $3\nmid k$ with conductor $2k-2$.

The simple parameterization with semigroup $\langle 3,k \rangle$ are classified in \cite{MP1}, 
and it was shown that (see  \cite[lemma 4,5,7,8 ]{MP1})

\begin{enumerate}
\item
$ (x(t),y(t))= (t^3, t^k +\sum_{i\geq k+1} a_it^i)$
is not simple if $k\geq 2p+9$. 
\item
If $p+9\leq k<2p+9$ then $ (x(t),y(t))= (t^3, t^k +t^l+\sum_{i\geq l+1} a_it^i)$,\\ $kl\equiv 2 \text{ mod } 3$, is not simple if $l\geq k+p$ and $k(k+p)\equiv 2 \text{ mod } 3$.
\end{enumerate}
 We will prove in this section that parameterizations with semigroup $\langle 3,k\rangle$ are unimodal in case of (1) or (2) and derive the normal form of Table \ref{tab.p>4,m<5}.\\

Let $l > k\geq 4$ such that $kl\equiv 2 \text{ mod } 3$. The gaps of the semigroup $\langle 3, k\rangle$ 
greater than $l$ are $\{l+3i \mid 1\leq i \leq \frac{2k-l}{3}-1\}$. \\

Using Theorem \ref{thm.Zar} and lemma 4 and lemma 7 from \cite{MP1} we have to prove the following lemma \ref{lem.m=3}.

\begin{lem}\label{lem.m=3} 
\begin{enumerate} 
\item Let $ k< l< 2k-9$ and $p \mid l-k$.
If $a\in K$, $a\neq 0$, then we have for $i>1$
$$(x(t),y(t))=(t^3,t^{k}+t^{l}+at^{l+3}+bt^{l+3i})\sim_\mA(t^3,t^{k}+t^{l}+at^{l+3}).$$
\item If $k \geq 2p+9$ then  $(x(t),y(t))= (t^3, t^k +\sum_{i\geq k+1} a_it^i)$ is unimodal (for any choice of $a_i\in K$).
\item If $p+9\leq k<2p+9$ then $ (x(t),y(t))= (t^3, t^k +t^l+\sum_{i\geq l+1} a_it^i)$ with $kl\equiv 2 \text{ mod } 3$, is unimodal if $l\geq k+p$ and $k(k+p)\equiv 2 \text{ mod } 3$.
\end{enumerate}
\end{lem}

\begin{proof}
(1) By Lemma 4 and Lemma 7 in \cite{MP1} one has to consider   only the case $p \mid l-k$. We give explicitly the $\mathcal A$-equivalence 
modulo $t^{l+3i+1}$ to remove the term $bt^{l+3i}$ and inductively the higher order terms in $y(t)$
($\varphi$ and $(H,L)$  as in (\ref{no.1})).
$$\varphi(t):=t+\alpha t^{3i+1} \text{ with } \alpha=\frac{b}{3a}.$$
$$H(x,y):=x(1+\alpha x^i)^3=x+3\alpha x^{i+1}+3\alpha^2 x^{2i+1} +\alpha^3 x^{3i+1}$$
$$L(x,y):=y(1+\alpha x^i)^{k}=\sum_{j=0}^{k}\binom {k} j \alpha^j x^{ij} y$$
With this definition we obviously have 
$$H((x(t),y(t))=H(t^3,t^{k}+t^{l}+at^{l+3}+bt^{l+3i})= t^3(1+\alpha t^{3i})^3=\varphi(t)^3.$$
Now modulo $ t^{l+3i+1}$ we obtain with $y'(t) = t^{k}+t^{l}+at^{l+3}$,
$$y'(\varphi(t))= \varphi^{k}+\varphi^{l}+a\varphi^{l+3}=$$
$$\sum_{j=0}^{k}\binom {k} j \alpha^j t^{k+3ij}+t^{l}+l\alpha
t^{l+3i}+at^{l+3}+a\alpha (l+3)t^{l+3i}.$$

On the other hand modulo $ t^{l+3i+1}$
$$L(t^3,t^{k}+t^{l}+at^{l+3}+bt^{l+3i})=$$
$$\sum_{j=0}^{k}\binom {k} j \alpha^j t^{k+3ij}+t^{l}+at^{l+3}+k\alpha t^{l+3i}
+(b+ka\alpha) t^{l+3i}.$$
We obtain as condition for the equivalence
\begin{center}
$b+ak\alpha=a(l+3)\alpha$ and $l\alpha=k\alpha$
\end{center}
In $K$ we have $l=k$  and since $\alpha =\frac{b}{3a}$, we obtain
modulo $ t^{l+3i+1}$
$$y'(\varphi(t))=L(t^3,t^{k}+t^{l}+at^{l+3}+bt^{l+3i}),$$
and hence $(x(\varphi(t))),y'(\varphi(t))=(H(x(t),y(t)),L(x(t),y(t))),$ proving (1).

(2) Repeated arguments as in (1)  show that $(x(t),y(t))$ 
can  be reduced to a 1-parametric $\Gamma$-full family. 
To prove unimodality, we use the following argument:
Since $(x(t), y(t))$  is not simple by \cite{MP1} if $k \geq 2p+9$, the family is modular (Definition \ref{def.modular}).  It follows from  \cite{MP1}  and (1) that parameterizations with semigroup  $\langle 3, k\rangle, k \geq 4$, are either simple or unimodal. Since parameterizations  with semigroup $\langle 2, k\rangle$ are simple, the unimodality follows now from Proposition \ref{prop.modular} (4).

(3) A similar argument as in (2) can be used for case (3).
\end{proof}
\bigskip


\subsection{Multiplicity 4}
\subsubsection{parameterizations with semigroup $\langle 4,6,2k+7\rangle$, $k\geq 1$.}
\text{}\\

We assume that the characteristic $p$ of $K$ is different from $2, 3$.\\
In \cite{MP1} it is proved that a parameterization with semigroup $\langle 4,6,2k+7\rangle$ is not simple if and only if $k\geq k_0 := \min \{k \text{ such that } p\mid 2k+7\}$, see (\ref{no.1}).
 In particular,  the parameterization $(t^4,t^6+t^{2k+1}+\sum_{j>2k+1}a_jt^j)$ is not simple if $p \mid 2k+7$, and simple with normal form $(t^4,t^6+t^{2k+1})$ 
 if  $p \nmid 2k+7$.
 \begin{prop}
parameterizations with semigroup $\langle 4,6,2k+7\rangle$ are unimodal 
  with normal form 
  $$(t^4,t^6+t^{2k+1}+at^{2k+3}),$$
 $a\in K$, if and only if $k\geq k_0$. 
 \end{prop}

Using Theorem \ref{thm.Zar}  and arguments as in the proof of  Lemma \ref{lem.m=3} (2) we have to prove the following lemma.

\begin{lem}
$(t^4,t^6+t^{2k+1}+at^{2k+3}+\sum_{j>2k+3}a_jt^j)\sim_\mA (t^4,t^6+t^{2k+1}+at^{2k+3})$
if $p\mid 2k+7$.
\end{lem}

\begin{proof}

The gaps of the semigroup $\langle 4,6,2k+7\rangle$ greater than $2k+3$ are $2k+5 \text{ and } 2k+9$. The conductor
of the semigroup is $2k+10$. Using Theorem \ref{thm.Zar} it is enough to consider the parameterization
$(t^4,t^6+t^{2k+1}+at^{2k+3}+bt^{2k+5})$. We explicitly give the $\mathcal A$-equivalence
$(t^4,t^6+t^{2k+1}+at^{2k+3}+bt^{2k+5})\sim_\mA (t^4,t^6+t^{2k+1}+at^{2k+3})$:
\begin{align*}
 \varphi(t)&=t+\alpha t^5 \text{ with } \alpha=\frac{b}{2k-5},\\
H(x,y)&=x(1+\alpha x)^4,\\
L(x,y)&=y(1+\alpha x)^6.
\end{align*}

Obviously $\varphi(t)^4=(t+\alpha t^5)^4=H(t^4,t^6+t^{2k+1}+at^{2k+3}+bt^{2k+5})$.\\
We obtain modulo $t^{2k+6}$
$$\varphi^6+\varphi^{2k+1}+a\varphi^{2k+3}=(t+\alpha t5)^6+t^{2k+1}+at^{2k+3}+(2k+1)\alpha t^{2k+5},$$
$$L(t^4,t^6+t^{2k+1}+at^{2k+3}+bt^{2k+5})=(t+\alpha t^5)^6+t^{2k+1}+at^{2k+3}+(b+6\alpha)t^{2k+5}.$$
This implies that modulo $t^{2k+6}$
$$\varphi^6+\varphi^{2k+1}+a\varphi^{2k+3}=L(t^4,t^6+t^{2k+1}+at^{2k+3}+bt^{2k+5}),$$
and proves the claim.
\end{proof}
\medskip

 \subsubsection{parameterizations with semigroup $\langle 4, 9\rangle$.}
\text{}\\

The semigroup  is
$\Gamma =\langle 4,9 \rangle =\{0,4,8,9,12,13,16,17,18,20,21,22,24,...\}$, with conductor $c=24$.
If $p=$ char$(K)\neq 2,3$ then, by Theorem \ref{thm.Zar}, the parameterization 
 is $\mA$-equivalent to
$$(x(t),y(t)) = (t^4, t^{9}+a_{10}t^{10} + a_{11}t^{11} + a_{14}t^{14} + a_{15}t^{15} +a_{19}t^{19}).$$
We will obtain the following normal forms:
 \begin{align*}
 & (t^4,t^9+t^{10}+at^{11}), \ a\in K, a\neq \frac{19}{18},\\
 & (t^4, t^{9}+(t^4, t^{9}+t^{10}+\frac{19}{18}t^{11}+bt^{15}), \ b\in K,\\
& (t^4,t^9+t^{11}),\\
&(t^4,t^9+t^{15}),\\
&(t^4,t^9+t^{19}),\\
&(t^4,t^9).
 \end{align*}
\begin{prop}\label{prop.4,9}
These are exactly the  unimodal parameterizations with semigroup $\Gamma = \langle 4, 9\rangle$ if $p \neq 2,3.$
  \end{prop}
 
  \begin{proof}
$\mathbf{Case\text{ } 1:}$ $a_{10}\neq 0$. 

Then, by replacing $t$ by $\alpha t$, dividing $x(t)$ by $\alpha ^4$, $y(t)$ by $\alpha ^9$ and setting $\alpha = \frac{1}{a_{10}}$
we may assume that $a_{10}=1$.

\begin{lem} \label{lem.a14}
If $a_{11}\neq 0, \frac{19}{18}$ then 
$$(x(t),y(t)) \sim_\mA  (t^4, t^{9}+t^{10}+at^{11}+bt^{19}).$$
\end{lem}
\begin{proof}
If $a_{14}\neq 0$ we will kill the term $a_{14}t^{14}$ in $y(t)$, leaving $x(t)$.
By substituting $t$ by $\varphi(t)=t-\frac{a_{14}}{9}t^6$ in $(x(t),y(t))$, we obtain ($\dots =$  terms of order $\geq 15$)
$$(\bar x(t),\bar y(t))=(t^4-\frac{4}{9}a_{14}t^9+\frac{2}{27}a_{14}^2t^{14}+ \dots, \ t^{9}+t^{10}+a_{11}t^{11}+ \dots)$$
Now we add $\frac{4}{9}a_{14}\bar y(t)$ to $\bar x(t)$ and obtain
$$(t^4+\alpha_{10}t^{10}+...+\alpha_{14}t^{14}+\dots, \ t^{9}+t^{10}+a_{11}t^{11}+\dots)$$
for suitable $\alpha_{10},...\alpha_{14}$. We substitute $t$ by $\varphi(t)=t-\frac{\alpha_{10}}{4}t^7+\dots$ to obtain
$$(t^4, t^{9}+t^{10}+a_{11}t^{11}+ \text{ terms of order } \geq 15 ),$$
This proves that we can assume $a_{14}=0$ and we do so in the following.\\

Now we  prove that we can reach $a_{15}=0$ by constructing an appropriate  $\mathcal A$-transformation using {\sc Singular}.\footnote{
\,We show the input of  {\sc Singular} in typewriter font and  give detailed explanations 
(only here, in the first place where we use {\sc Singular}). {\tt >} indicates output of {\sc Singular}.}
\begin{verbatim}
ring R=(0,a,b,c,d,a1,a2,a3,a4,a5,a6,a7,a8,a9,a10,b1,b2,u1,u2,u3,u4,
        u5,u6,v1,v2,v3,v4,v5,v6),(x,y,t),ds;
\end{verbatim}
The computation takes place in the ring {\tt R}, the localization of the polynomial ring $\Q(a,b,...,v_6)[x,y,t]$ with variables $x,y,t$ and 28 parameters
{\tt a,b,...,v6}. 
\begin{verbatim}
poly phi=t+a2*t2+a3*t3+a4*t4+a5*t5+a6*t6+a7*t7+a8*t8+a9*t9+a10*t10
         +b1*t11+b2*t12;
poly Hx=u1*x+u2*y+u3*x2+u4*xy+u5*x3+u6*x4;
poly Hy=v1*x+v2*y+v3*x2+v4*xy+v5*x3+v6*x4;
\end{verbatim}
The polynomial {\tt phi} ($\varphi$) is an automorphism of $K[[t]]$ and the pair {\tt (Hx,Hy)} ($H_x,H_y$) defines an automorphism of $K[[x,y]]$ if $u_1v_2 - u_2v_1 \neq 0$. The automorphisms depend on parameters  {\tt ai,bi,ui,vi} that have to be determined in order to ensure that $a_{15}=0$ is possible with such automorphisms. Here $K=\Q$. 
\begin{verbatim}
poly W=jet(phi^4-subst(Hx,x,t4,y,t9+t10+a*t11+b*t15),15);
poly X=jet(phi^9+phi^10+c*phi^11-subst(Hy,x,t4,y,t9+t10+a*t11+b*t15),15);
\end{verbatim}
We have $W= j^{15}(\varphi^4 - H_x(t^4,t^9+t^{10}+at^{11}+bt^{15}))$, a polynomial in $R$, and in the same way we get 
$X = j^{15}(\varphi^{9}+\varphi^{10}+c \varphi^{11} - H_y(t^4,t^9+t^{10}+at^{11}+bt^{15}))$. 
Note that we start with the parametrization $(x(t), y(t)) =(t^4,t^9+t^{10}+at^{11}+bt^{15})$, $a,b\neq0$ (higher order terms are irrelevant here).
Then we compare it with a parameterization $(x'(t),y'(t)) = (t^4,t^9+t^{10}+ct^{11})$ without term of degree 15. Here $a,b$ are given while $c, a_i, b_i,u_i,v_i$ are free parameters.
The condition $W=X=0$ means 
$$(x'(\varphi(t)), y'(\varphi(t))) = (H_x(x(t), y(t)),H_y(x(t), y(t))) \text{ mod } \langle t^{16}\rangle$$
or, equivalently, that $(x(t),y(t))$ is $\mathcal A$-equivalent to a parameterization without term of degree 15. Comparing coefficients of the above equation, we will get the required conditions  on the parameters  {\tt c,ai,bi,ui,vi} that guarantee $a_{15}=0$.

\begin{verbatim}
matrix M1=coef(W,t);
matrix M2=coef(X,t);
ideal I;
int ii;
for(ii=1;ii<=ncols(M1);ii++){I[size(I)+1]=M1[2,ii];}
for(ii=1;ii<=ncols(M2);ii++){I[size(I)+1]=M2[2,ii];}
\end{verbatim}
The matrices $M_1$ and $M_2$ collect the coefficients of $W$ and $X$ together with the corresponding monomial in $t$ and the ideal $I$ collects only the coefficients. 
\begin{verbatim}
I;      //We show only the relevant generators of I
I[9]=(4*a2^3+12*a2*a3+4*a4) //We show only the relevant generators of I
I[10]=(6*a2^2+4*a3)
I[11]=(4*a2)
I[12]=(-u1+1)
I[18]=(9*a2-v2+1)
I[19]=(-v2+1)
I[20]=(-v3)
I[21]=(-v1)
\end{verbatim}
This implies that in characteristic different from $2$ and $3$,  $v_1=0, v_3=0, v_2=1, u_1=1, a_2=0, a_3=0, a_4=0$.
We set {\tt v1=0, v3=0, v2=1, u1=1, a2=0, a3=0, a4=0 }
\begin{verbatim}
phi=t+a5*t^5+a6*t^6+a7*t^7+a8*t^8+a9*t^9+a10*t^10+b1*t^11+b2*t^12;
Hx =x+u2*y+u3*x^2+u4*x*y+u6*x^3+u7*x^4;
Hy =y+v4*x*y+v6*x^3+v7*x^4;
W=jet(phi^4-subst(Hx,x,t4,y,t9+t10+a*t11+b*t15),15);
X=jet(phi^9+phi^10+c*phi^11-subst(Hy,x,t4,y,t9+t10+a*t11+b*t15),15);
\end{verbatim}
and recompute $I$ as above.
\begin{verbatim}
I[7]=(4*a6-u2)      //the relevant generators
I[8]=(4*a5-u3)
I[10]=(10*a5+9*a6-v4)
I[11]=(9*a5-v4)
I[12]=(-v6)
I[13]=(-a+c)
\end{verbatim}
We set {\tt c=a, v6=0, v4=9*a5 u3=4*a5, u2=4*a6} and recompute $I$ as above.
\begin{verbatim}
phi=t+a5*t^5+a6*t^6+a7*t^7+a8*t^8+a9*t^9+a10*t^10+b1*t^11+b2*t^12;
Hx =x+4*a6*y+4*a5*x^2+u4*x*y+u6*x^3+u7*x^4;
Hy =y+9*a5*x*y+v7*x^4;
W=jet(phi^4-subst(Hx,x,t4,y,t9+t10+a*t11+b*t15),15);
X=jet(phi^9+phi^10+a*phi^11-subst(Hy,x,t4,y,t9+t10+a*t11+b*t15),15);
I;
I[5]=(-4*a*a6+4*a8)
I[6]=(-4*a6+4*a7)
I[8]=(a5+9*a6)
We set a7=a6, a8=a*a6, a5=-9*a6 and recompute I:

phi=t-9*a6*t^5+a6*t^6+a6*t^7+a*a6*t^8+a9*t^9+a10*t^10+b1*t^11+b2*t^12;
Hx =x+4*a6*y-36*a6*x^2+u4*x*y+u6*x^3+u7*x^4;
Hy =y-81*a6*x*y+v7*x^4;
W=jet(phi^4-subst(Hx,x,t4,y,t9+t10+a*t11+b*t15),15);
X=jet(phi^9+phi^10+a*phi^11-subst(Hy,x,t4,y,t9+t10+a*t11+b*t15),15);
I; 
I[3]=(-108*a6^2+4*a10-u4)
I[4]=(486*a6^2+4*a9-u6)
We set u6=486*a6^2+4*a9, u4=-108*a6^2+4*a10 and recompute:

phi=t-9*a6*t^5+a6*t^6+a6*t^7+a*a6*t^8+a9*t^9+a10*t^10+b1*t^11+b2*t^12;
Hx=x+4*a6*y-36*a6*x^2+(-108*a6^2+4*a10)*x*y+(486*a6^2+4*a9)*x^3+u7*x^4;
Hy=y-81*a6*x*y+v7*x^4;
W=jet(phi^4-subst(Hx,x,t4,y,t9+t10+a*t11+b*t15),15);
X=jet(phi^9+phi^10+a*phi^11-subst(Hy,x,t4,y,t9+t10+a*t11+b2*t15),15);
W;
(6*a6^2-4*a10+4*b1)*t^14+(-4*a*a10-4*b*a6+12*a6^2+4*b2)*t^15
X;
(-18*a*a6-b+19*a6)*t^15
\end{verbatim}

Setting  $b_1=-1/4*(6*a_6^2-4*a_{10}), b_2=-1/4(-4a*a_{10}-4b*a_6+12*a_6^2)$
we get:
$W=0$ and 
$X=(-b-18*c*a6+19*a6)*t^{15}$.
The free parameters are now $v_7, u_7, a_{10}, a_9$ and we can set them 0.
If we set $a_6=b/(19-18*c)$ we get $X=0$ and $W=0$, which proves the lemma.
\end{proof}

\begin{lem} 
If $a_{11}\neq 0$ then
\begin{enumerate}
  \item if $a_{11}\neq \frac{19}{18}$ then 
    $(t^4, t^{9}+\sum\limits_{\substack{i>9}} a_it^i)\sim_\mA (t^4, t^{9}+t^{10}+a_{11}t^{11})$,
  \item if $a_{11}= \frac{19}{18}$ then 
    $(t^4, t^{9}+\sum\limits_{\substack{i>9}} a_it^i)\sim _{\mathcal{A}} (t^4, t^{9}+t^{10}+\frac{19}{18}t^{11}+a_{15}t^{15})$.
\end{enumerate}
\end{lem}

\begin{proof}
The proof of the lemma uses again {\sc Singular} as in the previous lemma. 
 
(1) Since $a_{11}\neq\frac{19}{18}$ we can start with 
$(t^4, t^{9}+t^{10}+at^{11}+bt^{19})$ (Lemma \ref{lem.a14}).
\begin{verbatim}
ring R=(0,a,b),(x,y,t),ds;
poly phi=t-9*b/(-18*a+19)*t9+b/(-18*a+19)*t10+b/(-18*a+19)*t11
                +a*b/(-18*a+19)*t12;
poly Hx=x+4*b/(-18*a+19)*xy-36*b/(-18*a+19)*x3;
poly Hy=y-81*b/(-18*a+19)*x2*y;
poly W=jet(phi^4-subst(Hx,x,t4,y,t9+t10+a*t11+b*t19),19);
poly X=jet(phi^9+phi^10+a*phi^11-subst(Hy,x,t4,y,t9+t10+a*t11+b*t19),19);
W;
0
X;
0
\end{verbatim}
This shows that the above automorphisms kill the term  $bt^{19}$.\\

(2) Since $a_{11}= \frac{19}{18}$ we start with
$(t^4,t^9+t^{10}+19/18 t^{11}+at^{15}+bt^{19})$
\begin{verbatim}
ring R=(0,a,b,d),(x,y,t),ds;
poly phi=t-9*d*t^5+d*t^6+d*t^7+19*d/18*t^8+(-26244*d^2+3211*d)/2916*t^10
         +(-30618*d^2+3211*d)/2916*t^11+(52488*a*d-656100*d^2+61009*d)/52488*t^12
         +(-26244*a*d^2+3211*a*d+2916*b*d+10206*d^3+1026*d^2)/2916*t^16;
poly Hx=x+4*d*y-36*d*x^2+(-104976*d^2+3211*d)/729*x*y+(-26244*d^3+3515*d^2)/486*y^2
        +486*d^2*x^3+(52488*d^3-2869*d^2)/27*x^2*y+(-8748*d^3+56*d^2)/3*x^4;
poly Hy=y-81*d*x*y+(17496*d^2+133*d)/6*x^2*y+280*d/9*x^4;
poly W=jet(phi^4-subst(Hx,x,t4,y,t9+t10+19/18*t11+a*t15+b*t19),19);
poly X=jet(phi^9+phi^10+19/18*phi^11+a*(t15+15*(-9*d)*t19)
          -subst(Hy,x,t4,y,t9+t10+19/18*t11+a*t15+b*t19),19);         
W;
0
X;
(-39366*a*d-729*b+13122*d^2-1805*d)/729*t^19
\end{verbatim}

In $X$ we replaced $\varphi^{15}$ by $a(t^{15}+15(-9d)t^{19})$ which is $\varphi^{15}$ mod $t^{20}$ in order to speed up the {\sc Singular} computations.
Note that $1805=19^2\cdot 5, 13122=3^8\cdot 2, 729=3^6$, in characteristic different from $2$ and $3$, we can choose $d$ in such a way that the above term with $t^{19}$ vanishes.
\end{proof}

This finishes the case $a_{10}\neq0$. Now assume that $a_{10}=0$.
\medskip

$\mathbf{Case\text{ } 2:}$ $a_{10}=0$.
\begin{lem}

 If $a_{11}\neq0$ then $(t^4, t^{9}+\sum\limits_{\substack{i>9}} a_it^i)\sim _{\mathcal{A}} (t^4, t^{9}+t^{11})$.
\end{lem}

\begin{proof}
Using the $K^*$-action, we may assume that $a_{11}=1$. As before we may assume that $a_{12}=a_{13}=a_{14}=0$. We first prove that
$$(t^4, t^{9}+t^{11}+a_{15}t^{15}+\dots)\sim _{\mathcal{A}} (t^4, t^{9}+t^{11}+\overline{a}_{19}t^{19}\dots).$$
We use again {\sc Singular}.
\begin{verbatim}
ring R=(0,a),(x,y,t),ds;
poly phi=t+a/2*t5;
poly Hx=x+2*a*x2+3/2*a^2*x3;
poly Hy=y+9/2*a*xy;
poly W=jet(phi^4-subst(Hx,x,t4,y,t9+t10+t11+a*t15),15);
poly X=jet(phi^9+phi^11-subst(Hy,x,t4,y,t9+t11+a*t15),15);
W;
0
X;
0
\end{verbatim}
Now we assume that $a_{15}=a_{16}=a_{17}=a_{18}=0$. We have to prove that
$$(t^4, t^{9}+t^{11}+a_{19}t^{19}+\dots)\sim _{\mathcal{A}} (t^4, t^{9}+t^{11}+\overline{a}_{23}t^{23}\dots).$$
We use again {\sc Singular}.
\begin{verbatim}
ring R=(0,a),(x,y,t),ds;
poly phi=t+t^9+(a-2)/9*t^10+(a-2)/9*t^11+(a-2)/9*t^12;
poly Hx=x+(4*a-8)/9*x*y+4*x^3;
poly Hy=y+(a-2)*y^2+9*x^2*y;
poly W=jet(phi^4-subst(Hx,x,t4,y,t9+t10+t11+a*t19),19);
poly X=jet(phi^9+phi^11-subst(Hy,x,t4,y,t9+t11+a*t19),19);
W;
0
X;
0
\end{verbatim}
We can kill the term $\overline{a}_{23}t^{23}$ using Theorem  \ref{thm.Zar}.
\end{proof}
\medskip

$\mathbf{Case\text{ } 3:}$
 Now we can assume that $a_{10}=a_{11}=a_{12}=a_{13}=a_{14}=a_{16}=a_{17}=a_{18}=0$.

 \begin{lem}
 Assume that $a_{15}\neq 0$ then
  $(t^4, t^{9}+\sum\limits_{\substack{i>9}} a_it^i)\sim _{\mathcal{A}} (t^4, t^{9}+t^{15})$.
\end{lem}

\begin{proof}
Using the $K^*$-action, we may assume that $a_{15}=1$. We first prove that
$$(t^4, t^{9}+\sum\limits_{\substack{i>9}} a_it^i)\sim _{\mathcal{A}} (t^4, t^{9}+t^{15}+\overline{a}_{23}t^{23}+\dots).$$
We use again {\sc Singular}.
\begin{verbatim}
ring R=(0,a,d),(x,y,t),ds;
poly phi=t+d*t5;
poly Hx=x+4*d*x2+6*d^2*x3+(4*d^3)*x4;
poly Hy=y+9*d*xy+36*d^2*x2*y;
poly W=jet(phi^4-subst(Hx,x,t4,y,t9+t10+t15+a*t19),19);
poly X=jet(phi^9+phi^15-subst(Hy,x,t4,y,t9+t15+a*t19),19);
W;
0
X;
(-a+6*d)*t^19
\end{verbatim}
Choosing $d=\frac{a}{6}$ we  kill $t^{19}$ and with Theorem  \ref{thm.Zar} 
the term $\overline{a}_{23}t^{23}$.
\end{proof}
\medskip

$\mathbf{Case\text{ } 4:}$
Now we can assume that $a_{10}=a_{11}=a_{12}=a_{13}=a_{14}=a_{16}=a_{17}=a_{18}=a_{20}=a_{21}=a_{22}=0$.

\begin{lem}\label{lem.a19}
 Assume that $a_{19}\neq 0$ then
  $(t^4, t^{9}+\sum\limits_{\substack{i>9}} a_it^i)\sim _{\mathcal{A}} (t^4, t^{9}+t^{19})$.
\end{lem}

\begin{proof}
Using the $K^*$-action, we may assume that $a_{19}=1$. Using Theorem  \ref{thm.Zar}  we can kill the term $\overline{a}_{23}t^{23}$.
\end{proof}
So far we proved that parameterizations with semigroup $\Gamma =\langle 4,9\rangle$ form a either a $\Gamma$-full, 1-dimensional modular family or deform in such a family. To finish the proof of Proposition \ref{prop.4,9}, we recall  from  \cite{MP1}  that parameterizations with semigroup  smaller than $\langle 4, 9\rangle$ are simple. The claim follows now from Proposition \ref{prop.modular} (4).
\end{proof}
\medskip

\subsubsection{parameterizations with semigroup $\langle 4,10,2k+11\rangle$.}
\text{}\\

Consider the semigroup $\Gamma=\langle 4, 10, 2k+11 \rangle = \{0,4,8,10,12,...,2k+10,2k+11,2k+12,2k+14,2k+15,2k+16,2k+18,2k+19,...\}$,
$k\geq 5$, with conductor $2k+18$.
Then any parameterization with semigroup $\Gamma$  is $\mathcal{A}$-equivalent to 
$(t^4, t^{10}+t^{2k+1}+\sum_{ v>2k+1, v \notin \Gamma} b_vt^v)$ 
 by Theorem \ref{thm.Zar} and it is easy to see that every such parameterization has semigroup 
 $\Gamma$. We assume in this section that $p\neq 2$.

\begin{prop}\label{prop.4.10.2k+5}
Let $\psi(t) =(t^4,x(t))$ be a parameterization with $x(t) = t^{10}+t^{2k+1}+\sum_{v>2k+1, v \notin \Gamma} b_vt^v$, $k\geq5$.
\begin{enumerate}
  \item If $p \nmid 2k-9$ then $\psi$ is $\mathcal{A}$-equivalent to 
   $$\psi_{a}:=(t^4, t^{10}+t^{2k+1}+at^{2k+3}+ \sum_{v \geq 2k+7, v \notin \Gamma} b_vt^v).$$

   \item Let us write $\psi$ in the form
    $$\psi_{a,b}:=(t^4, t^{10}+t^{2k+1}+at^{2k+3} +bt^{2k+5}  + \sum_{v \geq 2k+7, v \notin \Gamma} b_vt^v).$$ 
If $p \mid 2k-9$ then $\psi_{a,b} \sim _{\mathcal{A}} \psi_{c,d}$  
if and only if $a=c$ and $b=d$ (i.e.  $\psi$ is at least bimodal).
\end{enumerate}
\end{prop}

\begin{proof}
We make the following Ansatz:
\begin{align*}
&\varphi(t):=t+\sum_{i\geq 2} a_it^i, a_i \in K,\\
&H_x:=u_1x+u_2y+u_3x^2+\dots,\\
&H_y:=v_1x+v_2y+v_3x^2+\dots,
\end{align*}
with $u_i, v_i \in K$ and ``...'' higher order terms in $x,y$.  

Assume that $(t^4,x(t)) \sim (t^4,y(t))$,  i.e. 
\begin{align}
\tag{$\dagger$}
\label{*}
\begin{split}
\varphi(t)^4 &=H_x(t^4,x(t)),\\
y(\varphi(t)) &= H_y(t^4,x(t)) 
\end{split}
\end{align}
with
$$x(t)=t^{10}+t^{2k+1}+at^{2k+3}+bt^{2k+5}+\sum\limits_{v>2k+5}c_vt^v,$$
$$y(t)=t^{10}+t^{2k+1}+ct^{2k+3}+dt^{2k+5}+\sum\limits_{v>2k+5}d_vt^v.$$

\begin{lem}\label{lem.4.10.2k+5} 
The equivalence of  $(t^4,x(t))$ and $(t^4,y(t))$ via {\em(\ref{*})} gives the following restrictions on the coefficients of $\varphi$ and $H_x, H_y$:\\
$a_3=0$, $a_l=0$ for $l$ even and $l\leq 2k-4$. $u_1=1, u_2=4a_7, u_3=4a_5, a_{2k-2}=a_7$, $v_1=v_3=0$, $v_2=1$, $v_4=10a_5$, $a=c$.
\end{lem}

\begin{proof} 
(1) We consider the first equation of (\ref{*}) and compare both sides.
The automorphism $\varphi$ satisfies
$$ \varphi(t)^4 =t^4+4t^3(a_2t^2+\dots)+6t^2(a_2t^2+\dots)^2+\dots,$$ giving the left side, while
$$H_x(t^4,x(t)) = H_x(t^4,t^{10}+\dots)=u_1t^4+u_2(t^{10}+\dots)+u_3t^8+\dots$$
gives the right side.

The term $4a_2t^5$ appears on the left side,  but no term with $t^5$ occurs on the right side, implying $a_2=0$. 
There is a term $a_3t^6$ on the left side but $t^6$ does not occur on the right side, hence $a_3=0$.

We prove by induction that $a_l=0$ for $l$ even and $l\leq 2k-4$.  We proved already $a_2=0$.
Now let $l \leq 2k-4$ be even and
assume that $a_j=0$ for all $j$ even and $j\leq l$. If $l=2k-4$ we are done. If $l<2k-4$, we obtain a term $4a_{l+2}t^{l+5}$ with $l+5 < 2k+1$, on the left side, while $t^{l+5}$ does not occur on the right side.

Comparing again $H_x(t^4,x(t))$ with $\varphi(t)^4$ we obtain that $u_1=1, u_2=4a_7$ and $u_3=4a_5$. In $\varphi(t)^4$, we have the term $4a_{2k-2}t^{2k+1}$, since $a_l=0$ for even $l$ and $l\leq 2k-4$.  Moreover, since $u_2=4a_7$, we obtain $a_{2k-2}=a_7$.

Now compare the left side of the second equation of (\ref{*}) 
$$y(\varphi(t))=\varphi^{10}+\varphi^{2k+1}+c\varphi^{2k+3}+...= t^{10}+10a_5t^{14}+ ...+ t^{2k+1} + ... + ct^{2k+3} + ...$$

with the right side 
$$H_y(t^4,x(t)) =v_1 t^4+v_2(t^{10}+ t^{2k+1} + at^{2k+3} + ...) + v_3t^8 + 
v_4t^4(t^{10}+ t^{2k+1} +...) + ... .$$
By comparing relevant terms of the both sides,  we obtain as before that $v_1=0, v_2=1, v_3=0, v_4=10a_5$ and $a=c$.
\end{proof}

(2) To finish the proof of Proposition \ref{prop.4.10.2k+5}, consider the coefficients of $t^{2k+5}$ on both sides  of the second equation of (\ref{*}).  The left side is
$$y(\varphi(t))=\varphi^{10}+\varphi^{2k+1}+c\varphi^{2k+3}+d\varphi^{2k+5}+...,$$
with
\begin{flalign*}
 \varphi &= t + a_5t^5+...,& \\
\varphi^{10}&= t^{10}+...+10a_{2k-4}t^{2k+5}+ ..., &\\
\varphi^{2k+1}& = t^{2k+1}+(2k+1)a_5t^{2k+5}+..., &\\ 
\varphi^{2k+3} &= t^{2k+3} + ..., &\\
t^{2k+5}&= t^{2k+5}+....
\end{flalign*}
The right side is
$$H_y(t^4,x(t))=t^{10}+t^{2k+1} + at^{2k+3} +bt^{2k+5} +...+10a_5t^4(t^{10} +t^{2k+1} +...)+...$$
The coefficient of $t^{2k+5}$ on the left side is
$10a_{2k-4}+(2k+1)a_5+d =(2k+1)a_5+d$, since $a_{2k-4}=0$. The coefficient on the right side is $b+10a_5$. The condition to eliminate $dt^{2k+5}$ is hence  $b-d-(2k-9)a_5=0$. 

If  $p \mid 2k-9$ then $b=d$ and (2) follows with Lemma \ref{lem.4.10.2k+5}.
If $p \nmid 2k-9$ then $bt^{2k+5}$  can be eliminated by choosing $a_5=d/ (2k-9)$ and (1) follows.
\end{proof}

In Proposition \ref{prop.4.10.2k+5} we eliminated the term $t^{2k+5}$ if
$p \nmid 2k-9$.  Now we want to eliminate further terms.

\begin{prop}  \label{prop.4.10.2k+7}
Let $p \nmid 2k-9$ and $\psi$ as in Proposition \ref{prop.4.10.2k+5}. Then
$\psi$ is $\mathcal{A}$-equivalent to 
    $$\psi_{ a,b}:=(t^4, t^{10}+t^{2k+1}+at^{2k+3} +bt^{2k+7}  + \sum_{v > 2k+7, v \notin \Gamma} b_vt^v).$$ 
\begin{enumerate} 
      \item If $p \nmid 2k+11$  then $\psi$
 is $\mathcal{A}$-equivalent to 
  $$\psi_a:=(t^4, t^{10}+t^{2k+1}+at^{2k+3}).$$
  \item 
If $p \mid 2k+11$ then\footnote{If $p \mid 2k+11$ and $p \mid 2k-9$ then $p=2$ or $p=5$.}
 $\psi_{a,b} \sim _{\mathcal{A}} \psi_{c,d}$  
if and only if $a=c$ and $b=d$ (i.e.  $\psi$ is at least bimodal).
\end{enumerate}
\end{prop}

\begin{proof}The proof is similar to the proof of Proposition \ref{prop.4.10.2k+5}. 
By Proposition \ref{prop.4.10.2k+5}(1) $\psi$  is $\mathcal{A}$-equivalent to     $\psi_{a,b}$.
We choose $\varphi$ and $H_x, H_y$ as in the proof of Proposition \ref{prop.4.10.2k+5} and
$$x(t)=t^{10}+t^{2k+1}+at^{2k+3}+bt^{2k+7}+\sum\limits_{v>2k+7}c_vt^v,$$
$$y(t)=t^{10}+t^{2k+1}+ct^{2k+3}+dt^{2k+7}+\sum\limits_{v>2k+7}d_vt^v.$$

In the following we always use the results of Lemma \ref{lem.4.10.2k+5}. Comparing the coefficients of $t^{2k+5}$ on both sides of the second equation of (\ref{*}), we obtain $a_5(2k+1)$ on the left side and $10a_5$ on the right side. 
This implies $a_5(2k-9)=0$, hence $a_5=0$. 

We want to compare now the coefficients of $t^{2k+7}$ in both sides of the second equation of (\ref{*}).
The left side is
$$y(\varphi(t))=\varphi^{10}+\varphi^{2k+1}+a\varphi^{2k+3}+d\varphi^{2k+7}+...,$$
with
\begin{flalign*}
 \varphi &= t + a_7t^7+...,& \\
\varphi^{10} &= t^{10}+...+10a_{2k-2}t^{2k+7}+ ..., \\
\varphi^{2k+1} &=t^{2k+1}+(2k+1)a_7t^{2k+7}+..., \\ 
\varphi^{2k+3}&=t^{2k+3} + ..., \\
\varphi^{2k+7}&= t^{2k+7}+....
\end{flalign*}
The right side is
$$H_y(t^4,x(t))=t^{10}+t^{2k+1} + at^{2k+3} +bt^{2k+7} +...+10a_5t^4(t^{10} +t^{2k+1} +at^{2k+3}+...)+...$$

Comparing the coefficients of $t^{2k+7}$, we obtain
$$10a_{2k-2}+(2k+1)a_7+b=d.$$
Since $a_{2k-2}=a_7$, this gives $(2k+11)a_7+b-d=0$.

If  $p \mid 2k+11$ then $b=d$ and (2) follows.
If $p \nmid 2k+11$ then $bt^{2k+7}$  can be eliminated by choosing $a_7=d/ (2k+11)$.
\medskip

It remains to prove (1), i.e., that $\psi_a(t)$ is $\mathcal{A}$-equivalent to $(t^4,x(t))$ with $\psi_a(t)=(t^4,y(t))$ and
$$ y(t)=t^{10}+t^{2k+1}+at^{2k+3}.$$
By what we proved so far and since $2k+9, 2k+13, 2k+17$ are the non-considered gaps of the semigroup $\Gamma$, we may assume that

$$x(t)=t^{10}+t^{2k+1}+at^{2k+3}+c_{2k+9}t^{2k+9}+c_{2k+13}t^{2k+13} +c_{2k+17}t^{2k+17}$$
and we want to eliminate the terms with $t^{2k+9}, t^{2k+13}$ and
$t^{2k+17}$.

Let us compare the coefficients of $t^{2k+9}$ in both sides of the second equation of (\ref{*}).
The left side is
$$y(\varphi(t))=\varphi^{10}+\varphi^{2k+1}+a\varphi^{2k+3}.$$
with
\begin{flalign*}
 \varphi &= t + a_7t^7+..., &\\
\varphi^{10} &= t^{10}+...+10a_{2k}t^{2k+9}+ ..., \\
\varphi^{2k+1} &=t^{2k+1}+(2k+1)a_7t^{2k+7}+(2k+1)a_{2k+9}t^{2k+9}+..., \\ 
\varphi^{2k+3} &=t^{2k+3} + ....
\end{flalign*}
The right side is
\begin{align*}
H_y(t^4,x(t))&=x(t)+v_5x(t)^2 + v_6t^{12}+ v_7t^8x(t)+...\\
&=t^{10}+t^{2k+1} + at^{2k+3} +c_{2k+9}t^{2k+9} +...+v_7t^{2k+9}+....
\end{align*}
Comparing the coefficients of $t^{2k+9}$, we obtain
$$10a_{2k}+(2k+1)a_{2k+9}=c_{2k+9}+v_7.$$
Choosing $v_7$ properly, we may assume that $c_{2k+9}=0$.
\medskip

Let us consider now the term with $t^{2k+13}$ in $x(t)$. 
In order to remove this term we get in a similar way the condition 
$$10a_{2k+4} + (2k+1)a_{13}+(2k+3)aa_{11} = c_{2k+13} + 2av_5.
$$ 

 We can choose $v_5$ properly to remove the term  $c_{2k+13}t^{2k+13}$.

Thus we may assume that
$$x(t)=t^{10}+t^{2k+1}+at^{2k+3}+c_{2k+17}t^{2k+17}$$
and we want to eliminate the term with  $t^{2k+17}$.
This is a consequence of  Theorem \ref{thm.Zar}.
This proves (1) and finishes the proof of Proposition \ref{prop.4.10.2k+7}.
\end{proof}

\begin{cor}\label{cor.4.10}

Let $p\neq 2,7$. Then every parameterization
$$\psi(t)=(t^4, t^{10}+t^{2k+1}+\sum_{v>2k+1}a_vt^v)$$
with $k\geq\frac{p+9}{2}$ if $5\leq p\leq 19$,  resp. $k\geq\frac{p-11}{2}$ if $p\geq 23$,
 is at least bimodal.
\end{cor}

\begin{proof} 
If $p\leq 19$ and $k<\frac{p+9}{2}$ then the condition $p\mid 2k+11$ implies that $p=5$ or $p=7$.
If $p=7$ and $k=5$ then $p\mid 2k+11=21$ this implies that for all $k\geq 5$ the parameterization $\psi(t)=(t^4, t^{10}+t^{2k+1}+\sum_{v>2k+1}a_vt^v)$ is at least bimodal (\ref{prop.4.10.2k+7}).

Let $k\geq\frac{p+9}{2}$. Consider the deformation $\psi_s(t)=\psi (t)+ (0,st^{p+10})$, $s\in K$.  Then $\psi_s, s\neq 0,$ is at least bimodal by Proposition \ref{prop.4.10.2k+5}, since $p+10=2k' +1$ and hence $p\mid 2k'-9$.
If $k\geq\frac{p-11}{2}$ we consider the deformation $\psi_s(t)=\psi (t)+ (0,st^{p-10})$.  Then $\psi_s, s\neq 0,$ is at least bimodal by Proposition \ref{prop.4.10.2k+7}, since $p-10=2k' +1$ and hence $p\mid 2k'+11$ and
$p\nmid 2k'-9$.

By semicontinuity of the modality (Theorem \ref{thm.semicontmod}) $\psi$ itself is at least bimodal.
\end{proof}

 \subsubsection{parameterizations with semigroup $\langle 4, 11\rangle$.}
\text{}
\medskip

A parameterization with semigroup  $\langle 4, 11\rangle$ is not unimodal in characteristic $p=11$ (Lemma \ref{lem6.3}).
Assume that $p\neq 11$.
We will obtain the following normal forms
\begin{align*}
&(t^4,t^{11}+t^{13}+at^{14}),\\
&(t^4,t^{11}+t^{14}+at^{17}) \text{ if } a\neq \frac{25}{22},\\
&(t^4,t^{11}+t^{14}+\frac{25}{22}t^{17}+at^{21}),\\
&(t^4,t^{11}+t^{17}),\\
&(t^4,t^{11}+t^{21}),\\
&(t^4,t^{11}+t^{25}),\\
&(t^4,t^{11}).
\end{align*}
The gaps of the semigroup $\langle 4, 11\rangle$
greater than $11$ are $13,14,17,18,21,25,29$.
Using Theorem \ref{thm.Zar} we have to consider the parameterization
\begin{center}
$((x(t),y(t))=(t^4,t^{11}+a_{13}t^{13}+a_{14}t^{14}+a_{17}t^{17}+a_{18}t^{18}+a_{21}t^{21}+a_{25}t^{25}).$
\end{center}
\medskip

$\mathbf{Case\text{ } 1: }$ $a_{13}\neq 0$\\
As in Case 1 before Lemma \ref{lem.a14} we use the $K^{*}$ action to obtain $a_{13}=1$. 

$\bullet$ We use {\sc Singular} to  reduced $(t^4,t^{11}+t^{13}+at^{14}+bt^{17})$  to $(t^4,t^{11}+t^{13}+ct^{14})$ (for a detailed explanation see the proof of Lemma \ref{lem.a14}): 
\begin{verbatim}
ring R=(0,a,b),(x,y,t),ds;
poly phi=t+(b)/2*t^5+(-b^3)/8*t^13;
poly Hx=x+(2*b)*x^2+(3*b^2)/2*x^3;
poly Hy=y+(11*b)/2*x*y;
poly W=jet(phi^4-subst(Hx,x,t4,y,t11+t13+a*t14+b*t17),17);
poly X=jet(phi^11+phi^13+a*phi^14-subst(Hy,x,t4,y,t11+t13+a*t14+b*t17),17);
W;
> 0
X;
> 0
\end{verbatim}

$\bullet$ Let the characteristic $p$ be different from 2 and 11. Then we reduce 
 $(t^4,t^{11}+t^{13}+at^{14}+bt^{18})$ to $(t^4,t^{11}+t^{13}+at^{14})$:
\begin{verbatim}
ring R=(0,a,b),(x,y,t),ds;
poly phi=t+(b)/11*t^8+(b)/11*t^10+(a*b)/11*t^11+(19*b^2)/242*t^15;
poly Hx=x+(4*b)/11*y;
poly Hy=y;
poly W=jet(phi^4-subst(Hx,x,t4,y,t11+t13+a*t14+b*t18),18);
poly X=jet(phi^11+phi^13+a*phi^14-subst(Hy,x,t4,y,t11+t13+a*t14+b*t18),18);
W;
> 0
X;
> 0
\end{verbatim}

$\bullet$ Now reduce $(t^4,t^{11}+t^{13}+at^{14}+bt^{21})$  to $(t^4,t{11}+t^{13}+at^{14})$: 
\begin{verbatim}
ring R=(0,a,b),(x,y,t),ds;
poly phi=t+(b)/2*t^9;
poly Hx=x+(2*b)*x^3+(3*b^2)/2*x^5;
poly Hy=y+(11*b)/2*x^2*y;
poly W=jet(phi^4-subst(Hx,x,t4,y,t11+t13+a*t14+b*t21),21);
poly X=jet(phi^11+phi^13+a*phi^14-subst(Hy,x,t4,y,t11+t13+a*t14+b*t21),21);
W;
> 0
X;
> 0
\end{verbatim}

$\bullet$ Reduce  $(t^4,t^{11}+t^{13}+at^{14}+bt^{25}) $ to $(t^4,t^{11}+t^{13}+at^{14})$ in characteristic different from 2.
\begin{verbatim}
ring R=(0,a,b),(x,y,t),ds;        
poly phi=t+1/2*b*t13;
poly Hx=x+2*b*x4;
poly Hy=y+11/2*b*x3y;
poly W=jet(phi^4-subst(Hx,x,t4,y,t11+t13+a*t14+b*t25),25);
poly X=jet(phi^11+phi^13+a*phi^14-subst(Hy,x,t4,y,t11+t13+a*t14+b*t25),25);
W;
>0
X;
>0
\end{verbatim}
\medskip

$\mathbf{Case\text{ } 2: }$ 
$a_{13}= 0$,  $a_{14}\neq0$.
\medskip

As above we may assume  $a_{14}=1$. 

$\bullet$
 Let $p\neq 3$, then $(t^4,t^{11}+t^{14}+at^{17}+bt^{18})$ will be reduced to $(t^4,t^{11}+t^{14}+at^{17})$.
\begin{verbatim}
ring R=(0,a,b),(x,y,t),ds;        
poly phi=t+(b)/3*t5;
poly Hx=x+(4b)/3*x2+(2b2)/3*x3+(4b3)/27*x4;
poly Hy=y+(11b)/3*xy;
poly W=jet(phi^4-subst(Hx,x,t4,y,t11+t14+a*t17+b*t18),18);
poly X=jet(phi^11+phi^14+a*t^17-subst(Hy,x,t4,y,t11+t14+a*t17+b*t18),18);
W;
> 0
X;
> 0
\end{verbatim}

$\bullet$ We reduce $(t^4,t^{11}+t^{14}+at^{17}+bt^{21})$ to $(t^4,t^{11}+t^{14}+at^{17})$.

$a \neq \frac{25}{22}$.
\begin{verbatim}
ring R=(0,a,b),(x,y,t),ds;        
poly phi=t+(11*b)/(66*a-75)*t^5+(-b)/(22*a-25)*t^8+(-b)/(22*a-25)*t^11
          +(-a*b)/(22*a-25)*t^14+(-3*b^2)/(968*a^2-2200*a+1250)*t^15 
          +(-22*a*b^2+22*b^2)/(484*a^2-1100*a+625)*t^18;
poly Hx=x+(-4*b)/(22*a-25)*y+(44*b)/(66*a-75)*x^2
         +(-44*b^2)/(484*a^2-1100*a+625)*x*y
         +(242*b^2)/(1452*a^2-3300*a+1875)*x^3 
         +(-484*b^3)/(31944*a^3-108900*a^2+123750*a-46875)*x^2*y
         +(5324*b^3)/(287496*a^3-980100*a^2+1113750*a-421875)*x^4
         +(14641*b^4)/(18974736*a^4-86248800*a^3
         +147015000*a^2-111375000*a+31640625)*x^5;
poly Hy=y+(121*b)/(66*a-75)*x*y+(6655*b^2)/(4356*a^2-9900*a+5625)*x^2*y;
poly W=jet(phi^4-subst(Hx,x,t4,y,t11+t14+a*t17+b*t21),21);
poly X=jet(phi^11+phi^14+a*phi^17-subst(Hy,x,t4,y,t11+t14+a*t17+b*t21),21);
W;
> 0
X;
> 0
\end{verbatim}

$\bullet$ Let $a =\frac{25}{22}$, characteristic different from 2,3,11.  
\begin{verbatim}
ring R=(0,a,b),(x,y,t),ds;        
poly phi=t+(b)/3*t12+(b)/3*t15+(25b)/66*t18+(ab)/3*t22;
poly Hx=x+(4b)/3*xy;
poly Hy=y+(11b)/3*y2;
poly W=jet(phi^4-subst(Hx,x,t4,y,t11+t14+25/22*t17+a*t^21+b*t25),25);
poly X=jet(phi^11+phi^14+25/22*phi^17+a*phi^21
       -subst(Hy,x,t4,y,t11+t14+25/22*t17+a*t^21+b*t25),25);
W;
> 0
X;
> 0
\end{verbatim}
$\bullet$ We reduce $(t^4,t^{11}+t^{14}+at^{17}+bt^{25})$ to $(t^4,t^{11}+t^{14}+at^{17}) $, $p\neq 3 $. 
\begin{verbatim}
ring R=(0,a,b),(x,y,t),ds;        
poly phi=t+(b)/3*t^12+(b)/3*t^15+(a*b)/3*t^18;
poly Hx=x+(4*b)/3*x*y;
poly Hy=y+(11*b)/3*y^2;
poly W=jet(phi^4-subst(Hx,x,t4,y,t11+t14+a*t17+b*t25),25);
poly X=jet(phi^11+phi^14+a*phi^17-subst(Hy,x,t4,y,t11+t14+a*t17+b*t25),25);
W;
> 0
X;
> 0
\end{verbatim}
\medskip

$\mathbf{Case\text{ } 3: }$ 
$a_{13}= 0$,  $a_{14}=0$, $a_{17}\neq 0$.

Using the $K^{*}$ action we obtain $a_{17}=1$. 

$\bullet$ We reduce 
$(t^4,t^{11}+t^{17}+bt^{18})$ to $(t^4,t^{11}+t^{17})$, $p\neq 2,11$.

\begin{verbatim}
ring R=(0,b),(x,y,t),ds;        
poly phi=t+1/11*b*t8+1/11*b*t14+19/242*b^2*t15;
poly Hx=x+4/11*b*y;
poly Hy=y;
poly W=jet(phi^4-subst(Hx,x,t4,y,t11+t17+b*t18),18);
poly X=jet(phi^11+phi^17-subst(Hy,x,t4,y,t11+t17+b*t18),18);
W;
> 0
X;
> 0
\end{verbatim}

$\bullet$ We reduce 
$(t^4,t^{11}+t^{17}+bt^{21})$ to $(t^4,t^{11}+t^{17})$, $p\neq 2,3$.
\begin{verbatim}
ring R=(0,b),(x,y,t),ds;        
poly phi=t+(b)/6*t^5+(-b^3)/216*t^13;
poly Hx=x+(2*b)/3*x^2+(b^2)/6*x^3+(-11*b^4)/1296*x^5;
poly Hy=y+(11*b)/6*x*y+(55*b^2)/36*x^2*y;
poly W=jet(phi^4-subst(Hx,x,t4,y,t11+t17+b*t21),21);
poly X=jet(phi^11+phi^17-subst(Hy,x,t4,y,t11+t17+b*t21),21);
W;
> 0
X;
> 0
\end{verbatim}

$\bullet$ We reduce 
$(t^4,t^{11}+t^{17}+bt^{25})$ to $(t^4,t^{11}+t^{17})$, $p\neq 2,3,11$.
\begin{verbatim}
ring R=(0,b),(x,y,t),ds;        
poly phi=t+(b)/6*t^9+(b)/66*t^13;
poly Hx=x+(2*b)/3*x^3+(2*b)/33*x^4+(b^2)/6*x^5+(b^2)/33*x^6;
poly Hy=y+(11*b)/6*x^2*y+(b)/6*x^3*y;
poly W=jet(phi^4-subst(Hx,x,t4,y,t11+t17+b*t25),25);
poly X=jet(phi^11+phi^17-subst(Hy,x,t4,y,t11+t17+b*t25),25);
W;
> 0
X;
> 0
\end{verbatim}
\medskip

$\mathbf{Case\text{ } 4: }$ $a_{13}= 0$,  $a_{14}=0$, $a_{17}= 0$, $a_{21}\neq 0$.

Using the $K^{*}$ action we obtain $a_{21}=1$.

$\bullet$ We reduce 
$(t^4,t^{11}+t^{21}+bt^{25})$ to $(t^4,t^{11}+t^{21})$, $p\neq 2,5$.
\begin{verbatim}
ring R=(0,b),(x,y,t),ds;        
poly phi=t+(b/10)*t5;
poly Hx=x+(2b)/5*x2+(3b2)/50*x3+(b3)/250*x4+(b4)/10000*x5;
poly Hy=y+(11b)/10*xy+(11b2)/20*x2y+(33b3)/200*x3y;
poly W=jet(phi^4-subst(Hx,x,t4,y,t11+t21+b*t25),25);
poly X=jet(phi^11+phi^21-subst(Hy,x,t4,y,t11+t21+b*t25),25);
W;
> 0
X;
> 0
\end{verbatim}
\medskip

$\mathbf{Case\text{ } 5: }$ 
If $a_{13}= 0$,  $a_{14}=0$, $a_{17}= 0$, $a_{21}= 0$, $a_{25}\neq 0$ then we obtain as before
$$(x(t),y(t))\sim_\mA (t^4,t^{11}+t^{25}).$$


 \subsection{Multiplicity 5}
 \subsubsection{Parameterizations with semigroup $\langle 5, 6\rangle$}
\text{}\\

The parameterization is not unimodal in characteristic $5$ by Lemma \ref{lem.5,6}(1). Assume that $p\neq 5$.
We obtain the following normal forms
\begin{align*}
&(t^5,t^6+t^8+at^9),\\
&(t^5,t^6+t^9),\\
&(t^5,t^6+t^{14}),\\
&(t^5,t^6).
\end{align*}
The gaps of the semigroup $\langle 5, 6\rangle$
greater than $6$ are $7,8,9,13,14,19$.

Using Theorem \ref{thm.Zar} we have to consider the parameterization
\begin{center}
$(t^5,t^{6}+a_{7}t^{7}+a_{8}t^{8}+a_{9}t^{9}+a_{13}t^{13}+a_{14}t^{14}+a_{19}t^{19}).$
\end{center}
\medskip

$\bullet$ We reduce first $(t^5,t^6+bt^7)$  to $(t^5,t^6)$ mod $t^8$.
\begin{verbatim}
ring R=(0,b),(x,y,t),ds;        
poly phi=t+(b)/6*t^2+(b^2)/9*t^3;
poly Hx=x+(5*b)/6*y;
poly Hy=y;
poly W=jet(phi^5-subst(Hx,x,t5,y,t6+b*t7),7);
poly X=jet(phi^6-subst(Hy,x,t5,y,t6+b*t7),7);
W;
> 0
X;
> 0
\end{verbatim}
We may assume that $a_7=0$ and consider 
\begin{center}
$(t^5,t^{6}+a_{8}t^{8}+a_{9}t^{9}+a_{13}t^{13}+a_{14}t^{14}+a_{19}t^{19}).$
\end{center}
The term $a_{19}t^{19}$ can be removed in all cases using Theorem \ref{thm.Zar}.
\medskip

$\mathbf{Case\text{ } 1: }$ $a_{8}\neq 0$\\
We use the $K^{*}$ action to obtain $a_{8}=1$.

$\bullet$ We show that $(t^5,t^6+t^8+at^9+bt^{14})$ is equivalent 
$(t^5,t^6+t^8+at^9)$. 
 \begin{verbatim}
 ring R=(0,a,b),(x,y,t),ds;
 poly phi=t+(b)/2*t7+(b)/2*t9+(ab)/2*t10;
 poly Hx=x+5/2*b*x*y;
 poly Hy=y+3*b*y^2;
 poly W=jet(phi^5-subst(Hx,x,t5,y,t6+t8+a*t9+b*t14),14);
 poly X=jet(phi^6+phi^8+a*phi^9-subst(Hy,x,t5,y,t6+t8+a*t9+b*t14),14);
W;
> 0
X;
> 0
\end{verbatim}

$\bullet$ We reduce $(t^5,t^6+t^8+at^9+bt^{13})$  to $(t^5,t^6+t^8+at^9)$.  
\begin{verbatim}
ring R=(0,a,b),(x,y,t),ds;        
poly phi=t-1/6*(-b)*t8;
poly Hx=x-5/6*(-b)*y2;
poly Hy=y;
poly W=jet(phi^5-subst(Hx,x,t5,y,t6+t8+a*t9+b*t13),13);
poly X=jet(phi^6+phi^8+a*phi^9-subst(Hy,x,t5,y,t6+t8+a*t9+b*t13),13);
W;
> 0
X;
> 0
\end{verbatim}
\medskip

$\mathbf{Case\text{ } 2: }$ $a_{8}= 0$,  $a_{9}\neq 0$\\
We use the $K^{*}$ action to obtain $a_{9}=1$.

$\bullet$ We reduce $(t^5,t^6+t^9+at^{13}+bt^{14})$ to $(t^5,t^6+t^9+at^{13})$  
\begin{verbatim}
ring R=(0,a,b),(x,y,t),ds;        
poly phi=t+1/3*b*t6;
poly Hx=x+5/3*b*x2;
poly Hy=y+2*b*xy;
poly W=jet(phi^5-subst(Hx,x,t5,y,t6+t9+a*t13+b*t14),14);
poly X=jet(phi^6+phi^9+a*phi^13-subst(Hy,x,t5,y,t6+t9+a*t13+b*t14),14);
W;
> 0
X;
> 0
\end{verbatim}
\medskip

$\mathbf{Case\text{ } 3: }$ $a_{8}= 0$,  $a_{9}=0$, $a_{13}=0$,  $a_{14}\neq0$.\\
We obtain as before the normal form $(t^5,t^6+t^{14})$.

 \subsubsection{parameterizations with semigroup $\langle 5, 7\rangle$.}
\text{}
\medskip

The parameterization is not unimodal in characteristic 5,7 by Lemma \ref{lem.5,6}(2).

We assume that $p\neq 5$ We obtain the following normal forms
\begin{align*}
&(t^5,t^7+t^8+at^{11}),\\
&(t^5,t^7+t^{11}+at^{13}),\\
&(t^5,t^7+t^{13}),\\
&(t^5,t^7+t^{18}),\\
&(t^5,t^7).
\end{align*}
The gaps of the semigroup $\langle 5, 7\rangle$
greater than $7$ are $8,9,11,13,16,18,23$.
By Theorem \ref{thm.Zar} we have to consider the parameterization
\begin{center}
$(t^5,t^{7}+a_{8}t^{8}+a_{9}t^{9}+a_{11}t^{11}+a_{13}t^{13}+a_{16}t^{16}+a_{18}t^{18}+a_{23}t^{23}).$
\end{center}
\bigskip

$\bullet$ We reduce $(t^5,t^7+at^8+bt9)$ to $(t^5,t^7+at^8)$ mod $t^9$. 
 \begin{verbatim}
ring R=(0,a,b),(x,y,t),ds;        
poly phi=t+(b)/7*t^3+(a*b)/7*t^4+(5*b^2)/49*t^5;
poly Hx=x+(5*b)/7*y;
poly Hy=y;
poly W=jet(phi^5-subst(Hx,x,t5,y,t7+a*t8+b*t9),9);
poly X=jet(phi^7+a*phi^8-subst(Hy,x,t5,y,t7+a*t8+b*t9),9);
 W;
> 0
X;
> 0
\end{verbatim}
We may assume that $a_9=0$ and consider\footnote{In characteristic $7$ this is not possible. Here we have to consider both cases, $a_9=0$ and  $a_9\neq 0$.
\begin{center}
$(t^5,t^{7}+a_{8}t^{8}+a_{11}t^{11}+a_{13}t^{13}+a_{16}t^{16}+a_{18}t^{18}+a_{23}t^{23}).$
\end{center}
}
\medskip

$\mathbf{Case\text{ } 1: }$ $a_{8}\neq 0$\\
We use the $K^{*}$ action to obtain $a_{8}=1$.

$\bullet$ Reduce $(t^5,t^7+t^8+at^{11}+bt^{13})$ to $(t^5,t^7+t^8+at^{11})$.
\begin{verbatim}
ring R=(0,a,b),(x,y,t),ds;        
poly phi=t+b*t6;
poly Hx=x+5*b*x2;
poly Hy=y+7*b*xy;
poly W=jet(phi^5-subst(Hx,x,t5,y,t7+t8+a*t11+b*t13),13);
poly X=jet(phi^7+phi^8+a*phi^11-subst(Hy,x,t5,y,t7+t8+a*t11+b*t13),13);
W;
> 0
X;
> 0
\end{verbatim}

$\bullet$ Reduce $(t^5,t^7+t^8+at^{11}+bt^{16})$  to $(t^5,t^7+t^8+at^{11})$.  
\begin{verbatim}
ring R=(0,a,b),(x,y,t),ds;        
poly phi=t+b*t8+b*t9+(a*b)*t12;
poly Hx=x+5*b*xy;
poly Hy=y+7*b*y2+b*x3;
poly W=jet(phi^5-subst(Hx,x,t5,y,t7+t8+a*t11+b*t16),16);
poly X=jet(phi^7+phi^8+a*phi^11-subst(Hy,x,t5,y,t7+t8+a*t11+b*t16),16);
W;
> 0
X;
> 0
\end{verbatim}

$\bullet$ Reduce $(t^5,t^7+t^8+at^{11}+bt^{18})$  to $(t^5,t^7+t^8+at^{11})$.  
\begin{verbatim}
ring R=(0,a,b),(x,y,t),ds;        
poly phi=t+(b)*t11;
poly Hx=x+(5*b)*x3;
poly Hy=y+(7*b)*x2y;
poly W=jet(phi^5-subst(Hx,x,t5,y,t7+t8+a*t11+b*t18),18);
poly X=jet(phi^7+phi^8+a*phi^11-subst(Hy,x,t5,y,t7+t8+a*t11+b*t18),18);
W;
> 0
X;
> 0
\end{verbatim}
\medskip

$\mathbf{Case\text{ } 2: }$ $a_{8}=0$,  $a_{11}\neq 0$\\
We use the $K^{*}$ action to obtain $a_{11}=1$.

$\bullet$ Reduce $(t^5,t^7+t^{11}+at^{13} +bt^{16})$  to $(t^5,t^7+t^{11}+at^{13})$.  
\begin{verbatim}
ring R=(0,a,b),(x,y,t),ds;        
poly phi=t+(1/4*b)*t6;
poly Hx=x+(5b)/4*x2+(5b2)/8*x3;
poly Hy=y+7/4*b*xy;
poly W=jet(phi^5-subst(Hx,x,t5,y,t7+t11+a*t13+b*t16),16);
poly X=jet(phi^7+phi^11+a*phi^13-subst(Hy,x,t5,y,t7+t11+a*t13+b*t16),16);
W;
> 0
X;
> 0
\end{verbatim}

$\bullet$ Show that $(t^5,t7+t^{11}+at^{13}+bt^{18})$ is equivalent $(t^5,t^7+t^{11}+at^{13})$.
\begin{verbatim}
 ring R=(0,a,b),(x,y,t),ds;
 poly phi =t+(b)/4*t8+(b)/4*t12+(ab)/4*t14;
 poly Hx=x+(5*b)/4*x*y;
 poly Hy=y+(7*b)/4*y^2;
 poly W=jet(phi^5-subst(Hx,x,t5,y,t7+t11+a*t13+b*t18),18);
 poly X=jet(phi^7+phi^11+a*phi^13-subst(Hy,x,t5,y,t7+t11+a*t13+b*t18),18);
W;
> 0
X;
> 0
\end{verbatim}
\medskip

$\mathbf{Case\text{ } 3: }$ $a_{8}=0$,  $a_{11}= 0$, $a_{13}\neq 0$\\
We use the $K^{*}$ action to obtain $a_{11}=1$.

$\bullet$ Reduce $(t^5,t^7+t^{13}+at^{16}+bt^{18})$ to $(t^5,t^7+t^{13}+at^{16})$.
\begin{verbatim}
ring R=(0,a,b),(x,y,t),ds;        
poly phi=t+(b)/6*t6;
poly Hx=x+(5b)/6*x2+(5b2)/18*x3;
poly Hy=y+(7b)/6*xy+(7b2)/12*x2y;
poly W=jet(phi^5-subst(Hx,x,t5,y,t7+t13+a*t16+b*t18),18);
poly X=jet(phi^7+phi^13+a*phi^16-subst(Hy,x,t5,y,t7+t13+a*t16+b*t18),18);
W;
> 0
X;
> 0
\end{verbatim}
\medskip

$\mathbf{Case\text{ } 4: }$ $a_{8}=0$,  $a_{11}= 0$, $a_{13}=0$, $a_{16}= 0$, $a_{18}\neq 0$\\
We obtain the following normal form
\begin{center}
$(t^5,t^7+t^{18}).$
\end{center}
\medskip

 \subsubsection{parameterizations with semigroup $\langle 5, 8\rangle$}
\text{}
\medskip

The parameterization is not unimodal in characteristic $5$ since parameterizations with semigroup $\langle 5, 8\rangle$ deform in parameterizations with semigroup $\langle 5, 7\rangle$, which is at least bimodal by Lemma \ref{lem.5,6}(1). We assume that $p\neq 5$. We obtain the following normal forms
\begin{align*}
&(t^5,t^8+t^9+at^{12}),\\
&(t^5,t^8+t^{12}+at^{14}),\\
&(t^5,t^8+t^{14}+at^{17}),\\
&(t^5,t^8+t^{17}),\\
&(t^5,t^8+t^{22}),\\
&(t^5,t^8).
\end{align*}
The gaps of the semigroup $\langle 5, 8\rangle$
greater than $8$ are $9,11,12,14,17,19,22,27$.

Using Theorem \ref{thm.Zar} we have to consider the parameterization
$$(t^5,t^{8}+a_{9}t^{9}+a_{11}t^{11}+a_{12}t^{12}+a_{14}t^{14}+a_{17}t^{17}+a_{19}t^{19}+a_{22}t^{22}+a_{27}t^{27}).$$
The term $a_{27}t^{27}$ can be removed in all cases by using Theorem \ref{thm.Zar}.\\

$\bullet$ Reduce $(t^5,t^8+at^9+bt^{11})$  to $(t^5,t^8+at^9)$ mod $t^{12} $.
\begin{verbatim}
ring R=(0,a,b),(x,y,t),ds;        
poly phi=t+(b)/8*t^4+(a*b)/8*t^5+(3*b^2)/32*t^7;
poly Hx=x+5/8*b*y;
poly Hy=y;
poly W=jet(phi^5-subst(Hx,x,t5,y,t8+a*t9+b*t11),11);
poly X=jet(phi^8+a*phi^9-subst(Hy,x,t5,y,t8+a*t9+b*t11),11);
W;
> 0
X;
> 0
\end{verbatim}

We may assume that $a_{11}=0$ and consider
$$(t^5,t^{8}+a_{9}t^{9}+a_{12}t^{12}+a_{14}t^{14}+a_{17}t^{17}+a_{19}t^{19}+a_{22}t^{22}+a_{27}t^{27}).$$
\medskip

$\mathbf{Case\text{ } 1: }$ $a_{9}\neq 0$\\
We use the $K^{*}$ action to obtain $a_{9}=1$.

$\bullet$ Reduce $(t^5,t^8+t^9+at^{12}+bt^{14})$  to $(t^5,t^8+t^9+at^{12})$.  
\begin{verbatim}
ring R=(0,a,b),(x,y,t),ds;        
poly phi=t+b*t6;
poly Hx=x+5*b*x2;
poly Hy=y+8*b*xy;
poly W=jet(phi^5-subst(Hx,x,t5,y,t8+t9+a*t12+b*t14),14);
poly X=jet(phi^8+phi^9+a*phi^12-subst(Hy,x,t5,y,t8+t9+a*t12+b*t14),14);
W;
> 0
X;
> 0
\end{verbatim}

$\bullet$ Reduce $(t^5,t^8+t^9+at^{12}+bt^{17})$  to $(t^5,t^8+t^9+at^{12})$.  
\begin{verbatim}
ring R=(0,a,b),(x,y,t),ds;        
poly phi=t+b*t9+b*t10+(a*b)*t13;
poly Hx=x+5*b*xy;
poly Hy=y+8*b*y2;
poly W=jet(phi^5-subst(Hx,x,t5,y,t8+t9+a*t12+b*t17),17);
poly X=jet(phi^8+phi^9+a*phi^12-subst(Hy,x,t5,y,t8+t9+a*t12+b*t17),17);
W;
> 0
X;
> 0
\end{verbatim}

$\bullet$ Reduce $(t^5,t^8+t^9+at^{12}+bt^{19})$  to $(t^5,t^8+t^9+at^{12})$.  
\begin{verbatim}
ring R=(0,a,b),(x,y,t),ds;        
poly phi=t+b*t11;
poly Hx=x+5*b*x3;
poly Hy=y+8*b*x2y;
poly W=jet(phi^5-subst(Hx,x,t5,y,t8+t9+a*t12+b*t19),19);
poly X=jet(phi^8+phi^9+a*phi^12-subst(Hy,x,t5,y,t8+t9+a*t12+b*t19),19);
W;
> 0
X;
> 0
\end{verbatim}

$\bullet$ Reduce $(t^5,t^8+t^9+at^{12}+bt^{22})$  to $(t^5,t^8+t^9+at^{12})$.  
\begin{verbatim}
ring R=(0,a,b),(x,y,t),ds;        
poly phi=t+(b)*t^14+(b)*t^15+(a*b)*t^18;
poly Hx=x+(5*b)*x^2*y;
poly Hy=y+(8*b)*x*y^2;
poly W=jet(phi^5-subst(Hx,x,t5,y,t8+t9+a*t12+b*t22),22);
poly X=jet(phi^8+phi^9+a*phi^12-subst(Hy,x,t5,y,t8+t9+a*t12+b*t22),22);
W;
> 0
X;
> 0
\end{verbatim}
\medskip

$\mathbf{Case\text{ } 2: }$ $a_{9}=0$,  $a_{12}\neq 0$\\
We use the $K^{*}$ action to obtain $a_{12}=1$.

$\bullet$ Reduce $(t^5,t^8+t^{12}+at^{14}+bt^{17})$  to $(t^5,t^8+t^{12}+at^{14})$. 
\begin{verbatim}
ring R=(0,a,b),(x,y,t),ds;        
poly phi=t+(b)/4*t^6;
poly Hx=x+(5*b)/4*x^2+(5*b^2)/8*x^3;
poly Hy=y+2*b*xy;
poly W=jet(phi^5-subst(Hx,x,t5,y,t8+t12+a*t14+b*t17),17);
poly X=jet(phi^8+phi^12+a*phi^14-subst(Hy,x,t5,y,t8+t12+a*t14+b*t17),17);
W;
> 0
X;
> 0
\end{verbatim}

$\bullet$ Reduce $(t^5,t^8+t^{12}+at^{14}+bt^{19})$  to $(t^5,t^8+t^{12}+at^{14})$. 
\begin{verbatim}
ring R=(0,a,b),(x,y,t),ds;        
poly phi=t+1/8*b*t12;
poly Hx=x+5/8*b*y2;
poly Hy=y;
poly W=jet(phi^5-subst(Hx,x,t5,y,t8+t12+a*t14+b*t19),19);
poly X=jet(phi^8+phi^12+a*phi^14-subst(Hy,x,t5,y,t8+t12+a*t14+b*t19),19);
 W;
> 0
X;
> 0
\end{verbatim}

$\bullet$ Reduce $(t^5,t^8+t^{12}+at^{14}+bt^{22})$  to $(t^5,t^8+t^{12}+at^{14})$. 
\begin{verbatim}
ring R=(0,a,b),(x,y,t),ds;        
poly phi=t+1/4*b*t11;
poly Hx=x+5/4*b*x3;
poly Hy=y+2*b*x2y;
poly W=jet(phi^5-subst(Hx,x,t5,y,t8+t12+a*t14+b*t22),22);
poly X=jet(phi^8+phi^12+a*phi^14-subst(Hy,x,t5,y,t8+t12+a*t14+b*t22),22);
W;
> 0
X;
> 0
\end{verbatim}
\medskip

$\mathbf{Case\text{ } 3: }$ $a_{9}=0$,  $a_{12}= 0$, $a_{14}\neq 0$\\
We use the $K^{*}$ action to obtain $a_{14}=1$.

$\bullet$ Reduce $(t^5,t^8+t^{14}+at^{17}+bt^{19})$  to $(t^5,t^8+t^{14}+at^{17})$. 
\begin{verbatim}
ring R=(0,a,b),(x,y,t),ds;        
poly phi=t+1/8*b*t12;
poly Hx=x+5/8*b*y2;
poly Hy=y;
poly W=jet(phi^5-subst(Hx,x,t5,y,t8+t14+a*t17+b*t19),19);
poly X=jet(phi^8+phi^14+a*phi^17-subst(Hy,x,t5,y,t8+t14+a*t17+b*t19),19);
W;
> 0
X;
> 0
\end{verbatim}

$\bullet$ Reduce $(t^5,t^8+t^{14}+at^{17}+bt^{22})$  to $(t^5,t^8+t^{14}+at^{17})$. 
\begin{verbatim}
ring R=(0,a,b),(x,y,t),ds;        
poly phi=t+(b)/6*t^9+(b)/6*t^15+(a*b)/6*t^18;
poly Hx=x+(5*b)/6*x*y+(5*b^2)/18*x*y^2;
poly Hy=y+(4*b)/3*y^2;
poly W=jet(phi^5-subst(Hx,x,t5,y,t8+t14+a*t17+b*t22),22);
poly X=jet(phi^8+phi^14+a*phi^17-subst(Hy,x,t5,y,t8+t14+a*t17+b*t22),22);
W;
> 0
X;
> 0
\end{verbatim}
\medskip

$\mathbf{Case\text{ } 4: }$ $a_{9}=0$,  $a_{12}= 0$, $a_{14}= 0$,$a_{17}\neq 0$\\
We use the $K^{*}$ action to obtain $a_{17}=1$.

$\bullet$ Reduce $(t^5,t^8+at^{17}+bt^{19})$  to $(t^5,t^8+t^{17})$. 

\begin{verbatim}
ring R=(0,a,b),(x,y,t),ds;        
poly phi=t+1/8*b*t12;
poly Hx=x+5/8*b*y2;
poly Hy=y;
poly W=jet(phi^5-subst(Hx,x,t5,y,t8+t17+b*t19),19);
poly X=jet(phi^8+phi^17-subst(Hy,x,t5,y,t8+t17+b*t19),19);
W;
> 0
X;
> 0
\end{verbatim}

We may assume that $a_{19}=0$ and consider
$$(t^5,t^{8}+t^{17}+a_{22}t^{22}+a_{27}t^{27}).$$

$\bullet$ Reduce $(t^5,t^8+t^{17}+bt^{22})$  to $(t^5,t^8+t^{17})$. 
\begin{verbatim}
ring R=(0,a,b),(x,y,t),ds;        
poly phi=t+1/9*b*t6;
poly Hx=x+(5*b)/9*x^2+(10*b^2)/81*x^3+(10*b^3)/729*x^4;
poly Hy=y+(8*b)/9*x*y+(28*b^2)/81*x^2*y;
poly W=jet(phi^5-subst(Hx,x,t5,y,t8+t17+b*t22),22);
poly X=jet(phi^8+phi^17-subst(Hy,x,t5,y,t8+t17+b*t22),22);
W;
> 0
X;
> 0
\end{verbatim}
\medskip

$\mathbf{Case\text{ } 5: }$ $a_{9}=0$,  $a_{12}= 0$, $a_{14}= 0$,$a_{17}= 0$, $a_{22}\neq 0$\\
As before we obtain as normal form
$$(t^5,t^8+t^{22}).$$

\section{Classification in characteristic $2,3$}
\label{sec.class3}
\subsection{Characteristic $2$}
In this section we assume characteristic $p=2$. 
\begin{lem} Let $p=2$. Then any unimodal parameterization is $\mA$-equivalent to one of the following normal forms with $a\in K$.
\begin{enumerate}
\item $(t^3,t^{k}+t^{l}+a\cdot t^{l+3})$, $8 \leq k< l \leq 2k-9$, $k\cdot l\equiv 2 \text{ mod } 3$,
\item $(t^3,t^{k})$,
\item $(t^4+t^6+a\cdot t^7,t^5)$,
\item $(t^4+t^6+a\cdot t^{11},t^5)$,
\item $(t^4+t^7,t^5)$,
\item $(t^4+ t^{11},t^5)$,
\item $(t^4,t^5)$,
\item $(t^4+t^5,t^6+a\cdot t^7)$,
\item $(t^4+t^5,t^6+a\cdot t^9)$,
\item $(t^4+t^5,t^6+a\cdot t^{11})$,
\item $(t^4+t^5,t^6+a\cdot t^{15})$,
\item $(t^4+t^5,t^6)$.
\end{enumerate}
\end{lem}

\begin{proof}
In \cite{MP1} it is proved that parameterizations with
semigroup $\Gamma< \langle 3,8\rangle$ are simple and parameterizations with semigroup $\langle 3,8\rangle$ are not simple. Lemma \ref{lem6.3} implies that parameterizations with semigroup $\Gamma\geq \langle 4,7 \rangle$ 
and the parameterization $(t^4+t^7+at^9,t^6+bt^7)$ are at least bimodal.\\
Since $2$ does not divide any denominator in the proof of lemma \ref{lem.m=3} we can prove (1) and (2)
as in characteristic different from $2$.\\
To prove (3) we consider the following computation:

\begin{verbatim}
ring R=(2,a,b),(x,y,t),ds;
poly phi=t+b/a*t5+b/a*t7;
poly Hx=x;
poly Hy=y+b/a*xy;
poly W=jet(phi^4+phi^6+a*phi^7+b*phi^11-subst(Hx,x,t4+t6+a*t7,y,t5),11);
poly X=jet(phi^5-subst(Hy,x,t4+t6+a*t7+b*t11,y,t5),11);
W;
> 0
X;
> 0
\end{verbatim}
The cases (4)-(7) are proved similarly.
To prove (8) we consider the following computations:

\begin{verbatim}
ring R=(2,a,e),(x,y,t),ds;
poly phi=t+a2*t2+e*a2*t3+(a*a2+a2^4)*t5;
poly Hx=x+a2*y+(a2^4)*x2;
poly Hy=y+(e*a2+a2^2)*x2;
poly W=jet(phi^4+phi^5-subst(Hx,x,t4+t5,y,t6+e*t7+a*t9),9);
poly X=jet(phi^6+e*phi^7-subst(Hy,x,t4+t5,y,t6+e*t7+a*t9),9);
W;
> 0
X;
> (a+e^2*a2+e*a2^2)*t^9  //will be 0 for an appropriate choice of e

ring R=(2,a,e),(x,y,t),ds;
poly phi=t+(e)*t2+(e2)*t3+(a+e5+e4)/(e2+e)*t4+(e4)*t5+(ae2+a+e4)/(e)*t7;
poly Hx=x+(e)*y+(a+e6+e4)/(e2+e)*x2+(a+e7+e4)/(e2+e)*xy;
poly Hy=y+(a)/(e+1)*xy;
poly W=jet(phi^4+phi^5-subst(Hx,x,t4+t5,y,t6+e*t7+a*t11),11);
poly X=jet(phi^6+e*phi^7-subst(Hy,x,t4+t5,y,t6+e*t7+a*t11),11);
W;
> 0
X;
> 0

ring R=(2,a,e),(x,y,t),ds;
poly phi=t+(e)*t^2+(e^2)*t^3+(e^3)*t^4+(e^4)*t^5+(e^3)*t^7
             +(a+e^10+e^9+e^8+e^6+e^5+e^4)/(e)*t^9
             +(a*e^3+a*e^2+a*e+a+e^13+e^11+e^10+e^4)/(e)*t^11;
poly Hx=x+(e)*y+(e^4+e^3)*x^2+(e^5+e^4+e^3)*x*y
             +(a+e^10+e^7+e^4)/(e)*y^2
             +(a+e^10+e^9+e^8+e^6+e^5+e^4)/(e)*x^3
             +(a*e^2+a+e^12+e^8+e^5+e^4)/(e)*x^2*y;
poly Hy=y+(e^7+e^6+e^4)*y^2+(e^7+e^4)*x^3+(e^9+e^7+e^6+e^4)*x^2*y;
poly W=jet(phi^4+phi^5-subst(Hx,x,t4+t5,y,t6+e*t7+a*t15),15);
poly X=jet(phi^6+e*phi^7-subst(Hy,x,t4+t5,y,t6+e*t7+a*t15),15);
W;
> 0
X;
> 0
\end{verbatim}
The cases (9)-(12) can be proved similarly.
\end{proof}
Normal forms of unimodal parameterizations are listed in the following table.
     \begin{table} [h!]  
  \caption{Unimodal Parametric Plane Curve Singularities,  $p=2$ (1 series, 11 sproadic singularities, $a\in K$)}  \label{tab.p=2}
        \begin{center}
    \begin{tabular}{| l | l |c|c|}
      \hline
       Normal form              & Semigroup    \\\hline
       
      \raisebox{-2 pt}{$(t^3,t^{k}+t^{l}+a\cdot t^{l+3}$), $8 \leq k< l \leq 2k-9$, $k\cdot l\equiv 2 \text{ mod } 3$}
                          & \raisebox{-2 pt}{$\langle 3,k \rangle$} \\
         $(t^3,t^{k})$  &  
        \\\hline
        
         \raisebox{-2 pt}{$(t^4+t^6+a\cdot t^7,t^5)$ }             &  \\
         \raisebox{-1 pt}{$(t^4+t^6+a\cdot t^{11},t^5)$}                    &$\langle 4,5 \rangle$ \\
        $(t^4+t^7,t^5)$  &  \\
        $(t^4+ t^{11},t^5)$  &  \\
        $(t^4,t^5)$   &  \\ \hline

         \raisebox{-2 pt}{$(t^4+t^5,t^6+a\cdot t^7)$}             &  \\
        $(t^4+t^5,t^6+a\cdot t^9)$                    &$\langle 4,6,13 \rangle$ \\
        $(t^4+t^5,t^6+a\cdot t^{11})$  &  \\
       $(t^4+t^5,t^6+a\cdot t^{15})$  &  \\
        $(t^4+t^5,t^6)$   &  \\ \hline

    \end{tabular}
   \end{center}
    \end{table}

\subsection{Characteristic $3$}

In this section we assume characteristic $p=3$. In \cite{MP1} it is proved that parameterizations with
with multiplicity $2$ are simple. parameterizations with semigroup $\langle 3,4\rangle$ are simple since the only
gap $\geq 3$ in the semigroup is $5$. parameterizations with semigroup $\Gamma\geq\langle 3,7\rangle$ are at least
bimodal (lemma \ref{lem6.3}) except parameterizations with semigroup $\langle 4,5\rangle$ which are simple (\cite{MP1}).

\begin{lem} Let $p=3$. Then any unimodal
parameterization has semigroup $\langle 3,5\rangle$  and  is $\mA$-equivalent to one of the following normal forms.
\begin{enumerate}
\item $(t^3+a\cdot t^4+t^7,t^5)$
\item $(t^3+t^4,t^5)$
\item $(t^3,t^5)$
\end{enumerate}
\end{lem}

\begin{proof}
The following computation shows that the parameterization $(t^3+a\cdot t^4+t^7,t^5)$ is not simple. Since the gaps of the semigroup $\geq 3$ are $4$ and $7$, the parameterization cannot have a higher modality than 1.

\begin{verbatim}
ring R=(3,a,c,a1,a2,a3,a4,a5,a6,a7,a8,a9,a10,b1,b2,b3,b4,b5,b6,
         b7,b8,b9,b10,a11,a14,a15,u1,u2,u3,u4,u5,u6,u7,u8,u9,u10,u11,
         u12,v1,v2,v3,v4,v5,v6,v7,v8,v9,v10,w1,w2,w3,w4,w5,w6,w7,w8,
         w9),(x,y,t),ds;
poly phi=t+a2*t2+a3*t3+a4*t4+a5*t5+a6*t6+a7*t7+a8*t8+a9*t9+a10*t10
        +b1*t11+b2*t12+b3*t13+b4*t14+b5*t15+b6*t16+b7*t17+b8*t18
        +b9*t19+b10*t20;
poly Hx=u1*x+u2*y+u3*x2+u4*xy+u5*y^2+u6*x3+u7*x2*y+u8*x4
              +u9*x5+u10*x6;
poly Hy=v1*x+v2*y+v3*x2+v4*xy+v5*y^2+v6*x3+v7*x2*y+v8*x4
              +v9*x5+v10*x6;
poly W=jet(phi^3+a*phi^4+phi^7-subst(Hx,x,t3+c*t4+t7,y,t5),11);
poly X=jet(phi^5-subst(Hy,x,t3+c*t4+t7,y,t5),11);
matrix M1=coef(W,t);
matrix M2=coef(X,t);
ideal I;
int ii;
for(ii=1;ii<=ncols(M1);ii++){I[size(I)+1]=M1[2,ii];}
for(ii=1;ii<=ncols(M2);ii++){I[size(I)+1]=M2[2,ii];}
ring S=3,(a1,a2,a3,a4,a5,a6,a7,a8,a9,a10,b1,b2,b3,b4,b5,b6,b7,b8,b9,
         b10,a11,a14,a15,u1,u2,u3,u4,u5,u6,u7,u8,u9,u10,u11,u12,v1,
    v2,v3,v4,v5,v6,v7,v8,v9,w1,w2,w3,w4,w5,w6,w7,w8,w9,a,b,c,d),lp;
ideal I=imap(R,I);
option(redSB);
ideal J=std(I);
J;
J[1]=a-c
J[2]=v3^3*c-v3^2*c^2-v4*c
J[3]=v2-1
J[4]=v1
...
\end{verbatim}
We get $a=c$ for appropriate choices of $a_i, b_i, u_i, v_i, w_i$ and hence that  $(t^3+a\cdot t^4+t^7,t^5)\sim_\mA  (t^3+c\cdot t^4+t^7,t^5)$ only if $a=c$.
\end{proof}
\medskip

Normal forms of unimodal parameterizations are listed in the following table.
     \begin{table} [h!]  
  \caption{Unimodal Parametric Plane Curve Singularities,  $p=3$    (3 sporadic singularities, $a\in K$)}
  \label{tab.p=3}
        \begin{center}
    \begin{tabular}{| l | l |c|c|}
      \hline
       Normal form              & Semigroup    \\\hline
         \raisebox{-2 pt}{$(t^3+a\cdot t^4+t^7,t^5)$}                &  \\
        $(t^3+t^4,t^5)$                   &$\langle 3,5 \rangle$ \\
        $(t^3,t^5)$  &  \\\hline

    \end{tabular}
   \end{center}
    \end{table}

\section{Non-unimodal parameterizations}
\label{sec.non}

In this section we will give a list of parameterizations which are not unimodal but in a sense at the border to unimodal parameterizations.
\begin{prop}\label{notuni}
Assume that the characteristic $p\neq 2$.
\begin{enumerate}
\item The parameterization $(t^4,t^{10}+t^{2k+1}+a\cdot t^{2k+3})$ is at least bimodal iff $k\geq \frac{p+9}{2}$.
\item A parameterization of multiplicity $4$ and semigroup $\Gamma\geq \langle 4,13\rangle$ is at least bimodal.
\item A parameterization with semigroup $\Gamma\geq \langle 5,9\rangle$ is tt least bimodal.
\end{enumerate}
\end{prop}
For proof the proposition we need the following lemma.

\begin{lem}\label{notunilem}
Assume that the characteristic $p\neq 2$.
\begin{enumerate}
\item The parameterization $(t^4,t^{10}+t^{p+10}+a\cdot t^{p+12})$ is at least bimodal. 
\item A parameterization with semigroup $\Gamma = \langle 4,13\rangle$ is at least bimodal.
\item A parameterization with semigroup $\Gamma = \langle 5,9\rangle$ is at least bimodal.
\item A parameterization with semigroup $\Gamma = \langle 6,7\rangle$ is at least bimodal.
\end{enumerate}
\end{lem}

\begin{proof} 
(1) is a consequence of Proposition \ref{prop.4.10.2k+5}.\\
(2) is a consequence of the following {\sc{Singular}}-computation.

\begin{verbatim}
int p=0;  \\p will be the characteristic in the following rings
ring R=(p,a,b,c,d,e,f,g,h,i,j,k,l,a2,a3,a4,a5,a6,a7,a8,a9,a10,b1,
            b2,b3,b4,b5,b6,b7,b8,b9,b10,a11,a14,a15,u1,u2,u3,u4,
            u5,u6,u7,u8,u9,u10,u11,u12,v1,v2,v3,v4,v5,v6,v7,v8,
            v9,w1,w2,w3,w4,w5,w6,w7,w8,w9),
            (x,y,t),ds;
poly phi=t+a2*t2+a3*t3+a4*t4+a5*t5+a6*t6+a7*t7+a8*t8+a9*t9+a10*t10
         +b1*t11+b2*t12+b3*t13+b4*t14+b5*t15+b6*t16+b7*t17+b8*t18
         +b9*t19+b10*t20;
poly Hx=u1*x+u2*y+u3*x2+u4*xy+u5*y^2+u6*x3+u7*x2*y+u8*x4+u9*x5;
poly Hy=v1*x+v2*y+v3*x2+v4*xy+v5*y^2+v6*x3+v7*x2*y+v8*x4+v9*x5;
poly W=jet(phi^4-subst(Hx,x,t4,y,t13+t14+a*t15+b*t19),19);
poly X=jet(ph^13+ph^14+c*ph^15+d*t^19
         -subst(Hy,x,t4,y,t13+t14+a*t15+b*t19),19);
matrix M1=coef(W,t);
matrix M2=coef(X,t);
ideal I;
int ii;
for(ii=1;ii<=ncols(M1);ii++){I[size(I)+1]=M1[2,ii];}
for(ii=1;ii<=ncols(M2);ii++){I[size(I)+1]=M2[2,ii];}
if(p==0)     
{
ring S=integer,(a2,a3,a4,a5,a6,a7,a8,a9,a10,b1,b2,b3,b4,b5,b6,b7,
                b8,b9,b10,a11,a14,a15,u1,u2,u3,u4,u5,u6,u7,u8,u9,
                u10,u11,u12,v1,v2,v3,v4,v5,v6,v7,v8,v9,w1,w2,w3,w4,
                w5,w6,w7,w8,w9,a,b,c,d,e,f,g,h,i,j,k,l),lp;
}
else
{    
ring S=p,(a2,a3,a4,a5,a6,a7,a8,a9,a10,b1,b2,b3,b4,b5,b6,b7,
                b8,b9,b10,a11,a14,a15,u1,u2,u3,u4,u5,u6,u7,u8,u9,
                u10,u11,u12,v1,v2,v3,v4,v5,v6,v7,v8,v9,w1,w2,w3,w4,
                w5,w6,w7,w8,w9,a,b,c,d,e,f,g,h,i,j,k,l),lp;
}            
ideal I=imap(R,I);
ideal J=std(I);
J;
//we display the first 4 entries of J
J[1]=8*b-8*d
J[2]=4*a-4*c
J[3]=4*v7
J[4]=2*v7^2-4*b+4*d
\end{verbatim}
This shows that $b=d$ and $a=c$ in characteristic different from 2 (and hence that the parameterization is at least bimodal).\\

From now on we do not give the whole {\sc{Singular}}-code, since it follows the same pattern as in (2). The rings are the same (except for the characteristic $p$), as well as the  polynomials Hx and Hy. From computation to computation we display only the changes, i.e., the characteristic, the polynomials W and X, and the result.\\

(3) is a consequence of the following {\sc{Singular}}-computation.
\begin{verbatim}
int p=0;
poly W=jet(phi^5-subst(Hx,x,t5,y,t9+t11+a*t12+b*t17) ,17);
poly X=jet(phi^9+phi^11+c*phi^12+d*t^17
                -subst(Hy,x,t5,y,t9+t11+a*t12+b*t17) ,17);
\end{verbatim}
\begin{verbatim}
J[1]=90*b-90*d
J[2]=5*a-5*c
J[3]=5*v6+40*b-40*d
J[4]=10*v4
\end{verbatim}
This shows that $b=d$ and $a=c$ in characteristic different from $2,3$ and $5$.\\

Now we consider the case in characteristic 5.
\begin{verbatim}
int p=5;
poly W=jet(phi^5+phi^6+f*phi^8+h*phi^11
                -subst(Hx,x,t5+t6+a*t8+c*t11,y,t9) ,14);
poly X=jet(phi^9-subst(Hy,x,t5+t6+a*t8+c*t11,y,t9) ,14);
\end{verbatim}

\begin{verbatim}
J[1]=c-h
J[2]=a-f
J[3]=v3
J[4]=v2-1
\end{verbatim}
This shows that $b=d$ and $a=c$ in characteristic $5$.\\

Consider now the case of characteristic $3$.

\begin{verbatim}
int p=3;
poly W=jet(phi^5-subst(Hx,x,t5,y,t9+t11+a*t12+b*t13) ,13);
poly X=jet(phi^9+phi^11+c*phi^12+d*t^13
                -subst(Hy,x,t5,y,t9+t11+a*t12+b*t13) ,13);
\end{verbatim}

\begin{verbatim}
J[1]=b-d
J[2]=a-c
J[3]=v3
J[4]=v2-1
\end{verbatim}
This shows that $b=d$ and $a=c$ in characteristic $3$.\\

(4) is a consequence of the following {\sc{Singular}}-computation.

\begin{verbatim}
int p=0;
poly W=jet(phi^6-subst(Hx,x,t6,y,t7+t9+a*t10+b*t11) ,11);
poly X=jet(phi^7+phi^9+c*phi^10+d*phi^11
                -subst(Hy,x,t6,y,t7+t9+a*t10+b*t11) ,11);  
\end{verbatim}

\begin{verbatim}
J[1]=42*b-42*d
J[2]=42*a-42*c
J[3]=6*a^2+6*b-6*c^2-6*d
J[4]=v2-1
\end{verbatim}
This shows that $b=d$ and $a=c$ in characteristic different from $2,3,7$.\\

Now we consider the case of characteristic 7.
\begin{verbatim}
int p=7;
poly W=jet(phi^6-subst(Hx,x,t6,y,t7+t8+a*t10+b*t11) ,11);
poly X=jet(phi^7+phi^8+c*phi^10+d*phi^11
                -subst(Hy,x,t6,y,t7+t8+a*t10+b*t11) ,11);
\end{verbatim}

\begin{verbatim}
J[1]=b-d
J[2]=a-c
J[3]=v2-1
J[4]=v1
\end{verbatim}
This shows that $b=d$ and $a=c$ in characteristic $7$.\\

Consider now the case of characteristic 3.
\begin{verbatim}
int p=3;
poly W=jet(phi^6+phi^8+c*phi^9+d*phi^10
                -subst(Hx,x,t6+t8+a*t9+b*t10,y,t7) ,10);
poly X=jet(phi^7-subst(Hy,x,t6+t8+a*t9+b*t10,y,t7) ,10);
\end{verbatim}

\begin{verbatim}
J[1]=b-d
J[2]=a-c
J[3]=v2-1
J[4]=v1
\end{verbatim}
This shows that $b=d$ and $a=c$ in characteristic $3$.
\end{proof}

\begin{proof} [Proof of Proposition \ref{notuni}]
Any parameterization from (1), (2) or (3) deforms to one of the  parameterizations listed in lemma \ref{notunilem}. The result follows now from the semicontinuity of the modality.
\end{proof}

 In the following we list some parameterizations which are not unimodal in small characteristic.
 \begin{lem}\label{lem6.3}
 \begin{enumerate}
  \item parameterizations with semigroup $\langle 3,7\rangle $ are at least bimodal in characteristic $3$.
 \item
 parameterizations with semigroup $\langle 4,11\rangle$ are not unimodal in characteristic
  $11$.

 \item parameterizations with semigroup $\langle 5,6\rangle$ are not unimodal in characteristic $2$.
 \item The parameterization $(t^4+t^7+a\cdot t^9,t^6+b\cdot t^7)$ is at least bimodal in characteristic $2$.
 \item parameterizations with semigroup $\Gamma \geq \langle 4,7\rangle$ are at least bimodal in characteristic $2$.
 \end{enumerate}
 \end{lem}
 
 \begin{proof}
 To prove (1) we consider the family $(t^3+\sum_{i\geq 4}a_it^i,t^7)$. Using the $K^{*}$-action we may assume that
$a_8=1$. Using theorem \ref{thm.Zar} we may consider the family $(t^3+at^4+bt^5+t^8,t^7)$. The result follows from the following computation:

\begin{verbatim}
int p=3;
poly W=jet(phi^3+a*phi^4+b*phi^5+phi^8
                -subst(Hx,x,t3+c*t4+d*t5+t8,y,t7),11);
poly X=jet(phi^7-subst(Hy,x,t3+c*t4+d*t5+t8,y,t7),11);
\end{verbatim}

\begin{verbatim}
J[1]=b-d
J[2]=a-c
J[3]=v4*c^2-v4*d-v6*c^3+v6*c*d
J[4]=v3
\end{verbatim}
This shows that $b=d$ and $a=c$ in characteristic $3$.\\

 To prove (2) we consider the following family of parameterizations $(t^4,t^{11}+t^{13}+a\cdot t^{14}+b\cdot t^{18})$.

\begin{verbatim}
int p=11;
poly W=jet(phi^4-subst(Hx,x,t4,y,t11+t13+a*t14+b*t18),18);
poly X=jet(phi^11+phi^13+c*phi^14+d*phi^18
                -subst(Hy,x,t4,y,t11+t13+a*t14+b*t18),18); 
\end{verbatim}

\begin{verbatim}
J[1]=b-d
J[2]=a-c
J[3]=v8
J[4]=v6
\end{verbatim}
This shows that $b=d$ and $a=c$ in characteristic $11$.\\

To prove (3) we consider the family $t^5,t^6+\sum_{i\geq 7}a_it^i$. Using the $K^{*}$-action we may assume that $a_8=1$. Now we use
following computations:

\begin{verbatim}
int p=2;
poly W=jet(phi^5-subst(Hx,x,t5,y,t6+c*t7+t8+d*t^9),13);
poly X=jet(phi^6+a*phi^7+phi^8+b*phi^9
                -subst(Hy,x,t5,y,t6+c*t7+t8+d*t^9),13);
\end{verbatim}

\begin{verbatim}
J[1]=b+d
J[2]=a+c
J[3]=v4
J[4]=v3^2+v3*c^4+v3*c^2+v3*c*d
\end{verbatim}
This shows that $b=d$ and $a=c$ in characteristic $2$.\\

To prove (4) we consider the following computations:

\begin{verbatim}
int p=2;
poly W=jet(phi^4+phi^7+a*t9-subst(Hx,x,t4+t7+d*t9,y,t6+e*t7),9);
poly X=jet(phi^6+b*phi^7-subst(Hy,x,t4+t7+d*t9,y,t6+e*t7),9);
\end{verbatim}

\begin{verbatim}
J[1]=b+e
J[2]=a*e+d*e
J[3]=v2+1
J[4]=v1
\end{verbatim}
This shows that $b=e$ and $a=d$ in characteristic $2$.\\

(5) is a consequence of (4).
 \end{proof}

\begin{cor}
parameterizations with semigroup $\Gamma \geq \langle 4,11\rangle$ are at least bimodal in characteristic 
$11$,
parameterizations with semigroup $\Gamma \geq \langle 3,7\rangle$ are at least bimodal in characteristic $3$ 
except parameterizations with semigroup $\langle 4,5\rangle$ which are simple.
parameterizations with semigroup $\Gamma \geq \langle 4,7\rangle$ are at least bimodal in characteristic $2$.
\end{cor}

\begin{lem}\label{lem.5,6}
\begin{enumerate}
\item
parameterizations with semigroup $\langle 5,6\rangle$ are at least bimodal in characteristic $3$ and $5$.
\item
parameterizations with semigroup $\langle 5,7\rangle$ are at least bimodal in characteristic $7$.
\end{enumerate}
\end{lem}

\begin{proof}
To prove (1) we consider the following computations.\\

First we consider the case of characteristic 3.

\begin{verbatim}
int p=3;
poly W=jet(phi^5-subst(Hx,x,t5,y,t^6+c*t7+t8+d*t9),9);
poly X=jet(phi^6+a*phi^7+phi^8+b*phi^9
                -subst(Hy,x,t5,y,t^6+c*t7+t8+d*t9),9);
\end{verbatim}

\begin{verbatim}
J[1]=b*c-c*d
J[2]=a-c
J[3]=v2-1
J[4]=v1
\end{verbatim}
This shows that $b=d$ and $a=c$ in characteristic $3$.\\

Now we consider the case of characteristic $5$.
\begin{verbatim}
poly W=jet(phi^5+phi^7+a*phi^8+b*phi^9
                -subst(Hx,x,t5+t7+c*t8+d*t9,y,t6) ,9);
poly X=jet(phi^6-subst(Hy,x,t5+t7+c*t8+d*t9,y,t6) ,9);
\end{verbatim}

\begin{verbatim}
J[1]=b-d
J[2]=a-c
J[3]=v2-1
J[4]=v1
\end{verbatim}
This shows that $b=d$ and $a=c$ in characteristic $3$.\\

To prove (2) we consider the following computation.
\begin{verbatim}
poly W=jet(phi^5-subst(Hx,x,t5,y,t7+a*t8+b*t9+t11),13);
poly X=jet(phi^7+e*phi^8+f*phi^9+phi^11
                -subst(Hy,x,t5,y,t7+a*t8+b*t9+t11),13);
\end{verbatim}
\begin{verbatim}
J[1]=b-f
J[2]=a-e
J[3]=v4*e^2+2*v4*f
J[4]=v3*e^2+2*v3*f
\end{verbatim}
This shows that $b=f$ and $a=e$ in characteristic $3$.
\end{proof}
\begin{cor}
\begin{enumerate}
\item
parameterizations with semigroup $\Gamma \geq\langle 5,6\rangle$ are at least bimodal in characteristic $5$.
\item
parameterizations with semigroup $\Gamma \geq\langle 5,7\rangle$ are at least bimodal in characteristic $7$.
\end{enumerate}
\end{cor}


\section{Results of Hefez and Hernandes for large characteristic} \label{sec.hefher}

In this section we first recall the theorem on $\mA$-normal forms  (Theorem \ref{Hefez}) of Hefez and Hernandes given in \cite{HHA},
who proved it for parametric plane branches over the complex numbers. The proof of Hefez-Hernandes uses the complete transversal theorem by Bruce, Kirk, and du Plessis for the action of complex Lie group on an affine space.

 We give a new proof of Theorem \ref{Hefez} without using the complete transversal theorem, which does not hold in positive characteristic.
Our proof follows the ideas of our classification  in characteristic $>0$ by constructing inductively the transformations of the group $ \mathcal A$ instead of working with the  Lie algebra. It applies to algebraically closed fields $K$ of characteristic 0 or of characteristic $p$, if $p$ is greater than the conductor $c(\Gamma)$ of the semigroup $\Gamma$ (Theorem  \ref{Hefez} and Corollary \ref{cor.p>c}).
A classification of the unimodal parametric plane branches in characteristic 0 by explicit lists of normal forms was first given in \cite{BM}, which we recall in Table \ref{tab.char0}. The classification in  \cite{BM} uses Theorem \ref{Hefez} and therefore holds also for $p>c(\Gamma)$.\\

Besides the use of the complete transversals Hefez and Hernandes
use the module of K\"ahler differentials and its semi-module of values to get more elimination criteria (as noticed already by O. Zariski). 

 Let $ (x(t), y(t)) \in K[[t]]^2$  be a parameterized plane branch with semigroup
$\Gamma $, conductor $c(\Gamma)$ and multiplicity $m(\Gamma)$. We set
$$\Omega:=\langle \frac{d(x(t))}{dt}, \frac{d(y(t))}{dt}\rangle_{K[[t]]} \subset K[[t]],$$
 $$\Gamma(\Omega):=\{\text{ord }_t(w)\mid w \in \Omega\} \text{ and  } $$ 
 $$\Lambda:=(\Gamma(\Omega)+1)\cup \{0\}.$$
$\Omega$ is isomorphic to the module of Kähler differentials of $K[[x(t),y(t)]]$ modulo torsion.
 $\Lambda$ is called the {\em semi-module (of values of differentials)} of the parameterization. We have $\Gamma\subseteq \Lambda$.
The number $\lambda:=$ min$\{s: s \in \Lambda\setminus \Gamma\}-m(\Gamma)$ is called the {\em Zariski number} of the parameterization. Note that $\lambda=\infty$ if $\Gamma=\Lambda$.  In Table  \ref{tab.char0} we give
$\lambda$ and the gaps $\{i \mid i>\lambda, i+ m(\Gamma)\not\in\Lambda\}$ of the semi-module $\Lambda$.\\

The following "Normal Forms Theorem" was proved by Hefez and Hernandes \cite[Theorem 2.1]{HHA} for $K=\C$. It improves Zariski's short parameterization to a "very short parameterization" by using the values of differentials.

\begin{thm}\label{Hefez}
Let $K$ be an algebraically closed field of characteristic $0$ and let $(x(t),y(t))$ be a parameterized plane branch. 
\begin{enumerate}
\item
There exist $\overline{x}(t)=t^n$, $\overline{y}(t)=t^m+t^{\lambda}+\sum_{i>\lambda, i+n \notin \Lambda}b_it^i$ such that
$$(x(t), y(t)) \sim _{\mathcal{A}} (\overline{x}(t), \overline{y}(t)).$$
\item If $(t^n,t^m+t^{\lambda}+\sum_{i>\lambda, i+n \notin \Lambda}a_it^i)\sim (t^n,t^m+t^{\lambda}+\sum_{i>\lambda, i+n \notin \Lambda}b_it^i)$ then there exists $r\in K^*$ such that $r^{\lambda-m}=1$ and $a_i=r^{i-m}b_i$ for all $i$.
\end{enumerate}
\end{thm}
The theorem holds if the characteristic of $K$ is $>c(\Gamma)$, see Corollary \ref{cor.p>c}.\\
  \begin{table} [h!]  
  \caption{Unimodal Parametric Plane Curve Singularities in 
  Characteristic $p$, $p=0$ \cite{BM}, \cite{MP} or $p>c(\Gamma)$ (1 series, 28 sporadic singularities, $a\in K$).
}  
  \label{tab.char0}
  
        \begin{center}
    \begin{tabular}{| l | l |c|c|}
      \hline
       Normal form              & Semigroup $\Gamma$   & $\lambda$  &Gaps of $\Lambda$ \\  \hline

       \raisebox{-2 pt}{$(t^4, t^9+t^{10}+at^{11})$, $a\neq \frac{19}{18}$}   & &10                      & 11 \\
        \raisebox{-1 pt}{$(t^4, t^9+t^{10}+\frac{19}{18}t^{11}+at^{15})$}       & &10              & 15 \\
       $(t^4, t^9+t^{11})$                    & $\langle 4,9 \rangle$ &11&-\\
       $(t^4, t^9+t^{15})$                     & &15&- \\
       $(t^4, t^9+t^{19})$                 & &19&- \\
       \raisebox{1 pt}{$(t^4, t^9)$}                 & & $\infty$ &-\\\hline

        \raisebox{-2 pt}{$(t^4, t^{10}+t^{2k+1}at^{2k+3})$}                     
        &  \raisebox{-1 pt}{$\langle 4,10,2k+11 \rangle$} &2k+1&2k+3
        \raisebox{-6 pt}{\text{}} \\
        \hline

        \raisebox{-2 pt}{$(t^4, t^{11}+t^{13}+at^{14})$}    &&13                   &14  \\
        $(t^4, t^{11}+t^{14}+at^{17})$, $a\neq \frac{25}{22}$     &&14                 & 17 \\
         \raisebox{-1 pt}{$(t^4, t^{11}+t^{14}+\frac{25}{22}t^{17}+at^{21})$}  &&14                  &21  \\
       $(t^4, t^{11}+t^{17})$                     &  $\langle 4,11 \rangle$ &17&-\\
       $(t^4, t^{11}+t^{21})$                 &  &21&-\\
       $(t^4, t^{11}+t^{25})$                        & &25&- \\
       $(t^4, t^{11})$                   & &$\infty$&- \\\hline

          \raisebox{-2 pt}{$(t^5, t^{6}+ t^{8}+at^{9})$}                & &8&9 \\
        $(t^5, t^6+t^9)$                    &$\langle 5,6 \rangle$ &9&-\\
        $(t^5, t^{6}+ t^{14})$  & &14&- \\
        $(t^5, t^{6})$   & &$\infty$&- \\\hline

        \raisebox{-2 pt}{$(t^5, t^{7}+t^{8}+at^{11})$ }               & &8&11 \\
        $(t^5, t^{7}+ t^{11}+at^{13})$                    &$\langle 5,7 \rangle$ &11&13\\
        $(t^5, t^{7}+ t^{13})$  & &13&- \\
        $(t^5, t^{7}+ t^{18})$  & &18&- \\
        $(t^5, t^{7})$   & &$\infty$&- \\\hline

         \raisebox{-2 pt}{$(t^5, t^{8}+t^{9}+at^{12})$}                & &9&12 \\
        $(t^5, t^{8}+ t^{12}+at^{14})$                    &$\langle 5,8 \rangle$ &12&14\\
        $(t^5, t^{8}+ t^{14}+at^{17})$  & &14&17 \\
        $(t^5, t^{8}+ t^{17})$  & &17&- \\
        $(t^5, t^{8}+ t^{22})$  & &22&- \\
        $(t^5, t^{8})$   & &$\infty$&- \\ \hline
    \end{tabular}
   \end{center}
    \end{table}

We give now a different proof of Theorem \ref{Hefez} by constructing inductively the $\mathcal A$-transformations $\varphi^{(k)}$ and $H^{(k)}$   up to order $k$. The following lemmas are formulated first for any $k$, using the exponential, but they are needed by Corollary  \ref{cor.p>c} only for $k\leq c(\Gamma)$ and are therefore valid if $p=0$ or if $p > c$. The following lemma is the basic inductive step.\\

Given $x(t)=(x_1(t),\dots, x_n(t)) \in K[[t]]^n$ a primitive parameterization
of a (not necessary plane) branch, $\varepsilon \in \langle t \rangle^2 K[[t]]$, $h=(h_1, \dots, h_n) \in \langle t \rangle^k K[[t]]^n$ and $g=(g_1, \dots, g_n) \in \langle x \rangle^2 K[[x]]^n$ with $x=(x_1,\dots,x_n)$. Let 
$$\delta=\varepsilon \frac{d}{dt}\text{ and  } \Delta=\sum g_i \frac{\partial}{\partial x_i}.$$
 Set
 $$H_i(x)=(\text{exp }\Delta)(x_i)=\sum^{\infty}_{v=0} \frac{1}{v!}\Delta^v(x_i)\text{ and } \varphi(t)=\text{exp }\delta(t)=\sum_{v=0}^\infty \frac{1}{v!}\delta^v(t)$$
 and
  $$H^{(k)}_i(x)=
  H_i \text{ mod } \langle x \rangle^k
  \text{ and } \varphi^{(k)}(t)=\text{exp }\delta(t)=\sum_{v=0}^k \frac{1}{v!}\delta^v(t).$$
    
\begin{lem}\label{lem1}
With these notions assume
 that $\delta(x(t))=\Delta(x)(x(t))+h \text{ mod } t^{k+1}$. Then, 
$$\varphi^{(k)}(x_i(t))=\varphi(x_i(t))=H_i(x(t)+h)=H^{(k)}_i(x(t)+h)\text{ mod }t^{k+1}.$$
\end{lem}

    \begin{proof}
      We know by assumption, with $x_i'(t) =  \frac{dx_i(t)}{dt}$, that\\
      $$\varepsilon \left(
    \begin{array}{c}
      x_1'(t) \\
                 . \\
                 . \\
                 . \\
                x_n'(t)
    \end{array}
  \right)
    = \left(
    \begin{array}{c}
      g_1(x(t))+h_1 \\
                 . \\
                 . \\
                 . \\
                g_n(x(t))+h_n
    \end{array}
  \right)      \,\,\,     mod \,\,\, t^{k+1}.$$
  This implies
  $$x(t)+\delta(x(t))=x(t)+h+\Delta(x)(x(t)+h)\,\,\, mod \,\,\, t^{k+1},$$
  since $h \in \langle t \rangle^k$, $g \in \langle x \rangle^2 K[[x]]$ and therefore
  $$\Delta(x)(x(t)+h)=\Delta(x)(x(t))\,\,\, mod \,\,\, t^{k+1}.$$
  Now we use the property
  $$\delta^2(x_i(t))=\delta(\Delta(x_i)(x(t))=\Delta^2(x_i)(x(t))\,\,\, mod \,\,\, t^{k+1},$$
  which we will prove in the next lemma.\\
  We obtain
  $$\sum\frac{1}{v!}\delta^v(x_i(t))=\sum\frac{1}{v!}\Delta^v(x_i)(x(t))\,\,\, mod \,\,\, t^{k+1},$$
  i.e.
  $$\varphi(x_i(t))=H_i(x(t)+h)\,\,\, mod \,\,\,t^{k+1}.$$

  Obviously we have
 $$\varphi^{(k)}(x_i(t))=\varphi(x_i(t))\text{ mod }t^{k+1} \text{ and } H_i(x(t)+h)=H^{(k)}_i(x(t)+h)\text{ mod }t^{k+1},$$
 since $\Delta^v(x_i)\in \langle x \rangle^{v+1}$ and $\delta^v(t)\in \langle t\rangle^{v+1}$.
 \end{proof}

\begin{lem}\label{lem2}
With the notations of Lemma \ref{lem1} the following hold:
\begin{enumerate}
  \item $\delta^2(x_i(t))=\delta(\Delta(x_i)(x(t))=\Delta^2(x_i)(x(t))\,\,\, mod \,\,\, t^{k+1}.$
  \item $\Delta^j(x_i)(x(t)+h)=\Delta^j(x_i)(x(t))\,\,\, mod \,\,\, t^{k+1}, j\geq 1.$
\end{enumerate}
\end{lem}

\begin{proof}
$\mathbf{(1)}$ $\Delta^2(x_i)=\Delta(g_i)=\sum_j g_j \frac{\partial g_i}{\partial x_j}$.
This implies that
 $$\Delta^2(x_i)(x(t)+h)=\sum g_j(x(t)+h) \frac{\partial g_i}{\partial x_j}(x(t)+h)$$
 $$                     =\sum g_j(x(t)) \frac{\partial g_i}{\partial x_j}(x(t))\,\,\, mod \,\,\, t^{k+1},$$
 since $h \in \langle t \rangle^k K[[x]]^n$. Moreover
 $$\delta(\Delta(x_i)(x(t))=\delta(g_i(x(t))=\sum_j \frac{\partial g_i}{\partial x_j}(x(t)) \cdot \delta((x_j(t))$$
 $$=\sum_j \frac{\partial g_i}{\partial x_j}(x(t)) (g_j(x(t))+h_j)\,\,\, mod \,\,\, t^{k+1}$$
 $$=\sum_j \frac{\partial g_i}{\partial x_j}(x(t)) (g_j(x(t)))\,\,\, mod \,\,\, t^{k+1},$$
 since $g_i \in \langle x \rangle^2 $ and $h_j \in \langle t \rangle^k $.\\
 This implies 
 $$\delta(\Delta(x_i)(x(t))=\Delta^2(x_i)(x(t))\,\,\, mod \,\,\, t^{k+1}.$$
 Now
 $$\delta^2(x_i(t))=\delta(g_i(x(t))+h_i)=\delta(g_i(x(t)))\,\,\, mod \,\,\, t^{k+1}$$
 $$=\delta(\Delta(x_i)(x(t))\,\,\, mod \,\,\, t^{k+1}$$
 since $\delta(h_i)=0 \,\,\, mod \,\,\, t^{k+1}$.\\
 $\mathbf{(2)}$ can be proved by using induction on $j$, since for $j=1$ we have
 $$g_i(x(t)+h)=g_i(x(t))\,\,\, mod \,\,\, t^{k+1}.$$
\end{proof}

\begin{exmp}
Let $x(t)=(t^5, t^9+t^{11})$, $\varepsilon=t^6$, $\Delta=5x_1^2\frac{\partial}{\partial x_1}+9x_1x_2\frac{\partial}{\partial x_2}$ then
$$\varepsilon \left(
    \begin{array}{c}
      x_1' \\\\

      x_2'
    \end{array}
  \right)
    = t^6\left(
    \begin{array}{c}
      5t^4 \\\\

     9t^8+11t^{10}
    \end{array}
  \right)
   = \left(
    \begin{array}{c}
      5t^{10} \\\\

     9t^{14}+11t^{16}
    \end{array}
  \right)  $$
  $$=\left(
    \begin{array}{c}
      5x_1^2(t) \\\\

      9x_1(t)x_2(t)+2t^{16}
    \end{array}
  \right)\,\,\,     mod \,\,\, t^{17},$$
  $\varphi(t)=t+t^6+3t^{11}+ \,\,\,     h.o.t. $, $H_1=x_1+5x_1^2+25x_1^3+ \,\,\,     h.o.t.$, $H_2=x_2+9x_1x_2+ \,\,\,     h.o.t.$
  and we have
  $$\varphi(x_1(t))=H_1(t^5, t^9+t^{11}+2t^{16}) \,\,\,     mod \,\,\, t^{17},$$
  $$\varphi(x_2(t))=H_2(t^5, t^9+t^{11}+2t^{16}) \,\,\,     mod \,\,\, t^{17}.$$
\end{exmp}

Now we weaken the assumptions of Lemma \ref{lem1} on $g$.

\begin{lem}\label{lem4}
We keep the assumptions of Lemma \ref{lem1} except just assuming $g\in \langle x \rangle K[[x]]^n$. Then 
$$\varphi(x_i(t)-h_i)=H_i(x(t))+\sum_l \frac{\partial g_i}{\partial x_l}(0) h_l\,\,\,     mod \,\,\, t^{k+1}.$$
and 
$$\varphi^{(k)}(x_i(t)-h_i)=H^{(k)}_i(x(t))+\sum_l \frac{\partial g_i}{\partial x_l}(0) h_l\,\,\,     mod \,\,\, t^{k+1}.$$
\end{lem}

\begin{proof}
Similar to the proof of Lemma \ref{lem2} we have
$$\delta(x(t))=\Delta(x)(x(t))+h\,\,\,     mod \,\,\, t^{k+1}.$$
This implies
$$\delta^2(x_i(t))=\delta \Delta(x_i)(x(t))$$
$$                =\sum_l \frac{\partial g_i}{\partial x_l} (x(t)) \delta x_l(t)$$
$$                =\sum_l \frac{\partial g_i}{\partial x_l} (x(t)) (g_l(x(t))+h_l)\,\,\,     mod \,\,\, t^{k+1}$$
$$                =\sum_l \frac{\partial g_i}{\partial x_l} (x(t)) g_l(x(t))+\sum_l \frac{\partial g_i}{\partial x_l}(0)h_l \,\,\,     mod \,\,\, t^{k+1}$$
since $h_l \in \langle t \rangle^k$ and $\frac{\partial g_i}{\partial x_l}(x(t))-\frac{\partial g_i}{\partial x_l}(0)\in \langle t \rangle$.\\
We obtain
$$\delta^2 x_i(t)=\Delta^2(x_i)(x(t))+\sum_l \frac{\partial g_i}{\partial x_l}(0)h_l \,\,\,     mod \,\,\, t^{k+1}.$$
Applying $\delta$ we obtain
$$\delta^j x_i(t)=\Delta^j(x_i)(x(t)) \,\,\,     mod \,\,\, t^{k+1}, j\geq 2.$$
This gives since $\delta(h_l)=0\,\,\,     mod \,\,\, t^{k+1}$ and $\delta x(t)=\Delta(x)(x(t))+h \,\,\,     mod \,\,\, t^{k+1}$ that
$(exp \,\, \delta)(x_i(t)-h_i)=(exp\,\,\Delta)(x_i)(x(t))+\sum_l \frac{\partial g_i}{\partial x_l}(0)h_l \,\,\,     mod \,\,\, t^{k+1}$.
The second property follows from the first by definition of $\varphi^{(k)}$ and $H^{(k)}$.
\end{proof}

\begin{exmp}
Let $x(t)=(t^5, t^8+t^9)$, $\delta=\frac{\alpha}{8}(t^4+t^5)dt$, $\Delta=\frac{5\alpha}{8}x_2\frac{\partial }{\partial x_1}$ then
$$\varphi(t)=(exp\,\,\delta)(t)=t+\frac{\alpha}{8}t^4+\frac{\alpha}{8}t^5+\frac{\alpha^2}{32}t^7+\,\,\,  h.o.t.$$
$H_1=x_1+\frac{5\alpha}{8}x_2$, $H_2=x_2$. We have
$$\varphi \left(
    \begin{array}{c}
      x_1(t) \\\\

      x_2(t)-\alpha t^{11}
    \end{array}
  \right)
    = exp\,\,\Delta\left(
    \begin{array}{c}
      x_1 \\\\

     x_2
    \end{array}
  \right)(x(t))
   + \left(
    \begin{array}{c}
      \frac{5\alpha^2}{16}t^{11} \\\\

     0
    \end{array}
  \right)  $$
  $$=\left(
    \begin{array}{c}
      H_1(x(t)) \\\\

     H_2(x(t))
    \end{array}
  \right)
   + \left(
    \begin{array}{c}
      \frac{5\alpha^2}{16}t^{11} \\\\

     0
    \end{array}
  \right) \,\,\,     mod \,\,\, t^{12}.$$
  Note that we can change $\varphi(t)=(exp\,\,\delta)(t)$ to $\psi(t)=t+\frac{\alpha}{8}t^5-\frac{\alpha^2}{32}t^7+\,\,\,     h.o.t.$ such that
  $$\psi\left(
    \begin{array}{c}
      x_1(t) \\\\

     x_2(t)-\alpha t^{11}
    \end{array}
  \right)
   =\left(
    \begin{array}{c}
      H_1(x(t)) \\\\

     H_2(x(t))
    \end{array}
  \right) \,\,\,     mod \,\,\, t^{12}.$$
  We will see later that this can be arranged in general.
\end{exmp}

We will now give a new proof of theorem \ref{Hefez} which follows the idea of our classification in characteristic $p>0$ (but without any computations that were needed in small characteristics). We assume now that $n=2$. The first part of the theorem is an immediate consequence of the following lemma. The second part can be proved as in \cite{HHA}.

\begin{lem}\label{lem6}
Given $x(t)=t^a$, $y(t)=t^b+t^{\lambda}+\sum_{i}a_it^i$. Assume that $k+a \in \Lambda$. Then there exists $a^{\prime}_i \in K$ such that $a^{\prime}_i = a_i$ for $i<k$ and $a^{\prime}_k =0$ and 
$$(x(t), y(t)) \sim _{\mathcal{A}} (t^a, t^b+t^{\lambda}+\sum_{i}a^{\prime}_it^i).$$
Moreover, let $\varphi$ and $H=(H_1,H_2)$ define the $\mathcal A$-equivalence above, i.e. 
$$\varphi(x(t),y(t))=H(t^a,t^b+t^{\lambda}+\sum_{i}a^{\prime}_it^i).$$
then
$$\varphi^{(k)}(x(t),y(t))\equiv H^{(k)}(t^a,t^b+t^{\lambda}+\sum_{i}a^{\prime}_it^i) \text{ mod } t^{k+1}.$$
\end{lem}

\begin{proof}
If $k+a \in \Lambda$ then there exist for $\alpha \in K$, $g_1, g_2 \in \langle x,y\rangle K[[x,y]]$ such that
$$a\cdot \alpha \cdot t^{k+a-1}=-g_2(x(t),y(t))\frac{d(x(t))}{dt}+g_1(x(t),y(t))\frac{d(y(t))}{dt}\,\,\,     mod \,\,\, t^{k+a}.\text{ (*) }$$
We obtain
$$a\cdot \alpha \cdot \frac{t^{k+a-1}}{x^{\prime}(t)}=\alpha t^k=-g_2(x(t),y(t))+\frac{g_1(x(t),y(t))}{x^{\prime}(t)}\cdot \frac{d(y(t))}{dt}\,\,\,     mod \,\,\, t^{k+a}.$$
We define $\delta=\frac{g_1(x(t),y(t))}{x^{\prime}(t)}\cdot \frac{d}{dt}$ and $\Delta=g_1\frac{\partial}{\partial x}+g_2\frac{\partial}{\partial y}$ and obtain
 $$\delta\left(
    \begin{array}{c}
      x(t) \\\\

     y(t)
    \end{array}
  \right)
   =\Delta\left(
    \begin{array}{c}
      x \\\\

     y
    \end{array}
  \right) (x(t), y(t))
  + \left(
    \begin{array}{c}
      0 \\\\

     \alpha t^k
    \end{array}
  \right) \,\,\,     mod \,\,\, t^{k+1}.$$
  Now we apply Lemma \ref{lem4} and obtain
  $$(exp\,\,\delta)\left(
    \begin{array}{c}
      x(t) \\\\

     y(t)-\alpha t^k
    \end{array}
  \right)
   =(exp\,\,\Delta)\left(
    \begin{array}{c}
      x \\\\

     y
    \end{array}
  \right) (x(t), y(t))
  + \left(
    \begin{array}{c}
     \frac{ \partial g_1}{\partial y}(0) \\\\

    \frac{ \partial g_2}{\partial y}(0)
    \end{array}
  \right)\alpha t^k \,\,\,     mod \,\,\, t^{k+1}.$$

  Now we have $\frac{ \partial g_2}{\partial y}(0)=0$
  and $\frac{g_1(x(t),y(t))}{x^{\prime}(t)}\in \langle t^2 \rangle$ (to obtain that $\delta(\langle t\rangle^l)\subseteq \langle t \rangle^{l+1}$).

  By using\,\footnote{\,This is a consequence  of $g_2 \in \langle x, y \rangle^2$ and $g_1\in Ky+\langle x, y \rangle^2$:
  Let $g_i=\alpha_ix+\beta_iy+ h.o.t.$ for $i=1,2$. Then since $a<b$ (*) implies that $\alpha_2=0$, $a\beta_2=b\alpha_1$ and $a\beta_2=\lambda\alpha_1$. The last two equations imply $\alpha_1=\beta_2=0$.
  Note that by definition of $\lambda$ always $2b>\lambda+a$.}
 $\frac{ \partial g_2}{\partial y}(0)=0$ and $\delta(\langle t\rangle^l)\subseteq \langle t \rangle^{l+1}$, we obtain with $\alpha = a_k$ that

  $$(exp\,\,\delta)\left(
    \begin{array}{c}
      x(t)-\alpha \frac{ \partial g_1}{\partial y}(0)t^k \\\\

     y(t)-a_k t^k
    \end{array}
  \right)
   =(exp\,\,\Delta)\left(
    \begin{array}{c}
      x \\\\

     y
    \end{array}
  \right) (x(t), y(t))\,\,\,     mod \,\,\, t^{k+1}$$
 i.e. $(x(t)-\alpha \frac{ \partial g_1}{\partial y}(0)t^k, y(t)-a_k t^k) \sim _{\mathcal{A}} (\overline{x}(t), \overline{y}(t))$ with $\overline{x}(t)=x(t)+$ terms of order $\geq k+1$ and $\overline{y}(t)=y(t)+$ terms of order $\geq k+1$. \\
 Now we perform the automorphism $\psi(t)=t+\frac{\alpha}{a}t^{k-a+1}+$ suitable terms of higher order smaller than $k+1$  and obtain
 $$(t^a, y(t)) \sim _{\mathcal{A}} (t^a, y(t)-a_kt^k\,\,\,  +   terms \,\, of \,\, order \,\, \geq k+1).$$
\end{proof}

\begin{cor}\label{cor.p>c}
Let $(x(t),y(t))=(t^a,t^b+\sum_{i>b}a_it^i)$ be a primitive parameterization with semigroup $\Gamma$, semi-module $\Lambda$, Zariski-number $\lambda$ and conductor $c=c(\Gamma)$.
Assume that the characteristic of $K$ is $p>c$. 
There exist $\overline{x}(t)=t^a$, $\overline{y}(t)=t^b+t^{\lambda}+\sum_{i>\lambda, i+a \notin \Lambda}b_it^i$ such that\,\footnote{\,The case $\lambda=\infty$ is included.}
$$(x(t), y(t)) \sim _{\mathcal{A}} (\overline{x}(t), \overline{y}(t)).$$
\end{cor}
\begin{proof}
We apply lemma \ref{lem6} to obtain 
$$(x(t), y(t)) \sim _{\mathcal{A}} (\overline{x}(t), \overline{y}(t)) \text{ mod } t^{c+1}.$$
Indeed, by Theorem \ref{thm.Zar}  and Remark  \ref{rm.Zar} (3)  the data $\varphi^{(k)}$ and $H^{(k)}$, which define the $\mathcal A$-equivalence, 
 are needed only for $k\leq c$ and then they are well defined, since $p>c$ by assumption. This finishes the proof.
\end{proof}
  \section*{Acknowledgments}
The research of the first author is supported by Higher Education Commission of Pakistan by the project no. 7495 /Punjab/NRPU/R$\&$D/HEC/2017.

 \section*{Conflict of Interests}
 On behalf of all authors, the corresponding author states that there is no conflict of interest. 

\end{document}